\documentclass[11pt,eps]{amsart}
\usepackage{amsmath}
\usepackage{amssymb,cmmib57,graphics}
\usepackage{graphicx}
\usepackage{times}
\usepackage{bm}
\usepackage{natbib}
\usepackage{threeparttable}

\newcommand{\RR}{\mathbb{R}}
\newcommand{\tr}{\text{tr}}

\newcommand{\calcommand}[2]{\newcommand{#1}{\mbox{$\mathcal #2$}}}
\calcommand{\calT}{T}
\calcommand{\calN}{N}
\calcommand{\calL}{L}
\newcommand{\E}{\text{E}}
\newcommand{\GP}{\text{GP}}
\newcommand{\diag}{\text{diag}}

\newcommand{\hcommand}[2]{\newcommand{#1}{\mbox{$\hat #2$}}}
\hcommand{\hgamma}{\gamma}
\hcommand{\hv}{v}

\newcommand{\hbcommand}[2]{\newcommand{#1}{\mbox{$\hat{\mathbf #2}$}}}
\hbcommand{\hbv}{v}

\newcommand{\tcommand}[2]{\newcommand{#1}{\mbox{$\tilde #2$}}}
\tcommand{\tbeta}{\beta}

\newcommand{\SP}{\text{SP}}
\newcommand{\SSR}{\mbox{SSR}}
\newcommand{\SSE}{\mbox{SSE}}
\newcommand{\bzero}{{\bf 0}}

\newcommand{\iidsim}{\stackrel{i.i.d.}{\sim}}

\newcommand{\cvgd}{\stackrel{d}{\rightarrow}}
\newcommand{\cvgp}{\stackrel{p}{\rightarrow}}

\newcommand{\dequ}{\stackrel{d}{=}}

\newcommand{\BFcommand}[2]{\def #1{{\mbox{\boldmath $#2$}}}}
\BFcommand{\balpha}{\alpha}
\BFcommand{\bbeta}{\beta}
\BFcommand{\bgamma}{\gamma}
\BFcommand{\bdelta}{\delta}
\BFcommand{\bepsilon}{\epsilon}
\BFcommand{\bvarepsilon}{\varepsilon}

\BFcommand{\bzeta}{\zeta}
\BFcommand{\bfeta}{\eta}
\BFcommand{\btheta}{\theta}
\BFcommand{\bvartheta}{\vartheta}
\BFcommand{\biota}{\iota}
\BFcommand{\bkapa}{\kapa}

\BFcommand{\blambda}{\lambda}
\BFcommand{\bmu}{\mu}
\BFcommand{\bnu}{\nu}
\BFcommand{\bxi}{\xi}
\BFcommand{\bpi}{\pi}
\BFcommand{\bvarpi}{\varpi}

\BFcommand{\brho}{\rho}
\BFcommand{\bvarrho}{\varrho}
\BFcommand{\bsigma}{\sigma}
\BFcommand{\bvarsigma}{\varsigma}
\BFcommand{\btau}{\tau}

\BFcommand{\bupsilon}{\upsilon}
\BFcommand{\bphi}{\phi}
\BFcommand{\bvarphi}{\varphi}
\BFcommand{\bchi}{\chi}
\BFcommand{\bpsi}{\psi}
\BFcommand{\bomega}{\omega}

\newcommand{\HBcommand}[2]{\newcommand{#1}{\hat{#2}}}
\HBcommand{\hbGamma}{\bGamma}

\newcommand{\bfcommand}[2]{\newcommand{#1}{{\mathbf #2}}}

\bfcommand{\ba}{a}
\bfcommand{\bb}{b}
\bfcommand{\bc}{c}
\bfcommand{\bd}{d}
\bfcommand{\bfe}{e}
\bfcommand{\bbf}{f}
\bfcommand{\bg}{g}
\bfcommand{\bh}{h}
\bfcommand{\bj}{j}
\bfcommand{\bk}{k}
\bfcommand{\bl}{l}

\bfcommand{\bn}{n}
\bfcommand{\bo}{o}
\bfcommand{\bp}{p}
\bfcommand{\bq}{q}
\bfcommand{\br}{r}

\bfcommand{\bs}{s}
\bfcommand{\bt}{t}
\bfcommand{\bu}{u}
\bfcommand{\bv}{v}
\bfcommand{\bw}{w}

\bfcommand{\bx}{x}
\bfcommand{\by}{y}
\bfcommand{\bz}{z}
\bfcommand{\bA}{A}
\bfcommand{\bB}{B}
\bfcommand{\bC}{C}
\bfcommand{\bD}{D}
\bfcommand{\bE}{E}
\bfcommand{\bF}{F}

\bfcommand{\bG}{G}
\bfcommand{\bH}{H}
\bfcommand{\bI}{I}
\bfcommand{\bJ}{J}
\bfcommand{\bK}{K}
\bfcommand{\bL}{L}

\bfcommand{\bM}{M}
\bfcommand{\bN}{N}
\bfcommand{\bO}{O}
\bfcommand{\bP}{P}
\bfcommand{\bQ}{Q}
\bfcommand{\bR}{R}

\bfcommand{\bS}{S}
\bfcommand{\bT}{T}
\bfcommand{\bU}{U}
\bfcommand{\bV}{V}
\bfcommand{\bW}{W}

\bfcommand{\bX}{X}
\bfcommand{\bY}{Y}
\bfcommand{\bZ}{Z}

\BFcommand{\bGamma}{\Gamma}
\BFcommand{\bDelta}{\Delta}
\BFcommand{\bTheta}{\Theta}
\BFcommand{\bLambda}{\Lambda}
\BFcommand{\bXi}{\Xi}
\BFcommand{\bPi}{\Pi}

\BFcommand{\bSigma}{\Sigma}
\BFcommand{\bUpsilon}{\Upsilon}
\BFcommand{\bPhi}{\Phi}
\BFcommand{\bPsi}{\Psi}
\BFcommand{\bOmega}{\Omega}

\theoremstyle{definition}
\newtheorem{theorem}{Theorem}[section]

\newtheorem{lemma}[theorem]{Lemma}
\newtheorem{thm}{Theorem}[section]

\theoremstyle{remark}

\newtheorem{proposition}[thm]{Proposition}

\title[A new Test for One-Way ANOVA with Functional Data]{A New Test for One-Way ANOVA with Functional Data and Application to Ischemic Heart Screening}

\author{Jin-Ting Zhang${}^\star$}
\address{Department of Statistics and Applied Probability, National University of Singapore, Singapore} \thanks{${}^\star$ Department of Statistics and Applied Probability, National University of Singapore, Singapore} \email{stazjt@nus.edu.sg} 

\author{Ming-Yen Cheng${}^\dagger$}
\address{Department of Mathematics, National Taiwan University, Taiwan} 
\thanks{${}^\dagger$ Department of Mathematics, National Taiwan University, Taiwan} 
\email{cheng@math.ntu.edu.tw}

\author{Chi-Jen Tseng${}^\ddagger$}
\address{Cardiovascular Research Institute, Fooyin University Hospital, Taiwan} 
\thanks{${}^\ddagger$ Cardiovascular Research Institute, Fooyin University Hospital, Taiwan} 
\email{chijen.tseng@gmail.com} 

\author{Hau-Tieng Wu${}^\Diamond$}
\address{Department of Mathematics, Stanford University, USA} 
\thanks{${}^\Diamond$ Department of Mathematics, Stanford University, USA} 
\email{hauwu@stanford.edu}

\begin{document}
\maketitle

\begin{abstract}
We propose and study a new global test, namely the $F_{\max}$-test,  for the one-way ANOVA problem in functional data analysis. The test statistic is taken as the maximum value of the usual pointwise $F$-test statistics over the interval the functional responses are observed. A nonparametric bootstrap method is employed to approximate the null distribution of the test statistic and to obtain an estimated critical value for the test. The asymptotic random expression of the test statistic is derived and the asymptotic power is studied. In particular, under mild conditions, the $F_{\max}$-test asymptotically has the correct level and is root-$n$ consistent in detecting local alternatives.  Via some simulation studies, it is found that in terms of both level accuracy and power,  the $F_{\max}$-test outperforms the Globalized Pointwise F (GPF) test of \cite{Zhang_Liang:2013} when the functional data are highly or moderately correlated, and its performance is comparable with the latter otherwise. An application to an ischemic heart real dataset suggests that, after proper manipulation, resting electrocardiogram (ECG) signals can be used as an effective tool in clinical ischemic heart screening, without the need of further stress tests as in the current standard procedure.
\end{abstract}

{\it keywords}  $F$-type test; $L^2$-norm based test; local power; myocardial ischemia; nonparametric bootstrap, pointwise $F$-test.

\section{Introduction}
Functional data are getting increasingly common in a wide scope of scientific fields ranging from climatology to medicine, from seismology to chemometrics, and so on.  In addition, as data collection technology evolves, nowadays we can conveniently collect and record a huge amount of functional data.  Compared with traditional data which consist of point observations, functional data might contain more detailed information about the underlying system. On the other hand, new challenges arise in the endeavor to extract any meaningful information hidden in the often massive functional data at hand. In the past two decades, functional data analysis has emerged as an important area and significant progress has been achieved  \citep{Ramsay_Silverman:2005,Ferraty:2011,Zhang:2013}.  In this paper we concentrate on the one-way ANOVA problem with functional responses, which is a fundamental problem in the inference but has received much less attention in functional data analysis compared to other problems such as regression. Besides, functional one-way ANOVA testing is closely related to classification, clustering and image analysis \citep{Xue_Titterington:2011,Slaets_Claeskens_Hubert:2012,Samworth:2012,Hall_Xia_Xue:2013}.

A particular medical example is about ischemic heart which is treated in greater detail in Section \ref{section:ischemia:screening}: we ask if we can distinguish the group of normal subjects from the group of subjects with ischemic heart by reading their resting electrocardiogram (ECG) signals, which are clearly functional data. Traditional, based on the physiological facts, physicians ``reduce the dimension'' of an ECG signal by focusing on a specific subinterval, and then distinguish between the two groups based on this reduced information. Intuitively, we lose information about the system by doing so, and a direct consequence is that we are not able to distinguish between the group of normal subjects from the group of stable ischemic heart  with the reduced data. Thus, in clinics, further stress tests are needed to assist the diagnosis.  However, stress tests cannot be applied to patients who are vulnerable to acute attacks during the testing, which is associated with the stress itself. An immediate question we may ask is whether we can avoid the stress. In particular, we ask if there is a better way to differentiate between the normal and ischemic heart groups if we process the resting ECG signals properly during the diagnosis. 
In other words, we like to test equality of the mean functions of the two groups of manipulated ECG functions, and this is a special case of the one-way ANOVA hypothesis testing problem for functional data.

Let $y_{i1}(t), y_{i2}(t), \cdots, y_{in_i}(t)$, where $i=1,2,\cdots, k$, denote $k$ groups of random functions defined over a given finite interval $\calT=[a,b]$ and $n_i\in \mathbb{N}$ is the number of cases in the $i$-th group. Let  $\SP(\mu, \gamma)$ denote a stochastic process with mean function $\mu(t), t\in\calT$, and covariance function $\gamma(s,t), s,t\in\calT$. Assuming that
\begin{equation}\label{ksamp.sec1}
y_{i1}, y_{i2},\cdots,y_{in_i}\iidsim \SP(\mu_i, \gamma), \; i=1,2,\cdots, k,
\end{equation} 
it is often interesting to test the equality of the $k$ mean functions:
\begin{equation}\label{anova.sec1}
H_0: \mu_1 (t)= \mu_2 (t) = \cdots = \mu_k (t), \; \forall t\in\calT,
\end{equation} 
against the usual alternative that at least two of the mean functions are not equal. The above problem is known as the $k$-sample testing problem or the one-way ANOVA problem for functional data. For the above one-way ANOVA problem, the $k$ mean functions are often decomposed as $\mu_i(t)=\mu_0(t)+\alpha_i(t), i=1,2,\cdots, k,$ where  $\mu_0(t)$ is the grand mean function and $\alpha_i(t), i=1,2,\cdots, k$, are  the $k$  main-effect functions so that the null hypothesis (\ref{anova.sec1}) is often  equivalently written as the problem of testing if the main-effect functions are all zero, i.e., $ H_{0} : \alpha_{1}(t) = \alpha_{2}(t) = \cdots = \alpha_{k}(t)= 0,\forall t\in\calT. $

There are several existing methods for testing the one-way ANOVA problem (\ref{anova.sec1}). The $L^2$-norm based test and the $F$-type test for linear models with functional responses 
\citep{Faraway:1997,Shen_Faraway:2004,Zhang_Chen:2007,Zhang:2011}
may be adopted for this purpose; see \cite{Zhang_Liang:2013} for more details.
\cite{Cuevas_Febrero_Fraiman:2004} proposed and studied  an $L^2$-norm based test directly for testing (\ref{anova.sec1}). They derived the limit random expression of their test statistic and suggested to approximate the null distribution by a parametric bootstrap method via re-sampling from the Gaussian process involved in the limit random expression under the null hypothesis.
In addition, the pointwise $F$-test  was proposed by \cite{Ramsay_Silverman:2005}, naturally extending the classical $F$-test to the context of functional data analysis; see also \cite{Zhang:2013}. It's test statistic is defined as
\begin{equation}\label{pointF.sec1}
F_n(t)=\frac{\SSR_n(t)/(k-1)} {\SSE_n(t)/(n-k)},
\end{equation} 
where and throughout,  $n=\sum_{i=1}^k n_i$ denotes the total sample size,
\begin{equation}\label{SSRSSE.sec1}
\SSR_n(t)=\sum_{i=1}^k n_i[{\bar{y}}_{i.}(t)-{\bar{y}}_{..}(t)]^2\,\, \mbox{ and }\,\,
\SSE_n(t)=\sum_{i=1}^k\sum_{j=1}^{n_i}[y_{ij}(t)-{\bar{y}}_{i.}(t)]^2
\end{equation} 
denote the  pointwise between-subject and within-subject variations respectively, and
\begin{equation}\label{ybar.sec1}
{\bar{y}}_{..}(t)=\frac{1}{n}\sum_{i=1}^k\sum_{j=1}^{n_i}y_{ij}(t)\,\,  \mbox{ and } \,\,
{\bar{y}}_{i.}(t)=\frac{1}{n_i}\sum_{j=1}^{n_i}y_{ij}(t),\, i=1,\ldots,k,
\end{equation} 
are the sample grand mean function and the sample group mean functions respectively. Then, the null hypothesis (\ref{anova.sec1}) is rejected as long as the pointwise null hypothsis $H_{0t}: \mu_1 (t)= \mu_2 (t) = \cdots = \mu_k (t)$ is rejected at any point $t\in\calT$.

There are a few advantages of using the above pointwise $F$-test.  Denote by $F_{p,q}$ the $F$-distribution with $p$ and $q$ degrees of freedom.
When the functional data (\ref{ksamp.sec1}) are realizations of Gaussian processes, under the null hypothesis (\ref{anova.sec1})   $F_n(t)$ follows the $F_{k-1,n-k}$ distribution for any given $t\in\calT$.  Therefore,  for  any pre-determined significance level $\alpha$, we can test the null hypothesis (\ref{anova.sec1})  at all of the points in $\calT$  using the same critical value $F_{k-1, n-k}(\alpha)$  where $F_{k-1,n-k}(\alpha)$ denotes the upper $100\alpha$ percentile of the $F_{k-1,n-1}$ distribution. On the other hand, the pointwise $F$-test has some limitations too. For example,  
for a given significance level, it is not guaranteed that alternative hypothesis in the one-way ANOVA
problem 
is overall significant even when the pointwise $F$-test is significant for each point in $\calT$ at the same significance level.  To overcome this difficulty, \cite{Cox_Lee:2008} proposed to correct the pointwise P-values of the pointwise $F$-tests via incorporating some  Bonferroni-type multiple comparison idea. The resulting pointwise $F$-test, however,  becomes very complicated since the
corrected pointwise P-values are obtained via intensive bootstrapping.  Therefore, some simple method for summarizing the pointwise $F$-tests is desirable. To this end, \cite{Zhang_Liang:2013} proposed and studied a so-called  Globalized Pointwise F-test, abbreviated as GPF test, whose test statistic is taken as the integral of the pointwise $F$-test statistic over $\calT$:
\begin{equation}\label{Tn.sec1}
T_n=\int_{\calT} F_n(t) dt.
\end{equation} 
Via  some simulation studies, \cite{Zhang_Liang:2013} showed that in terms of size-controlling and power,  the GPF test is in general comparable with  the $L^2$-norm based test and the $F$-type test adopted for the one-way ANOVA problem (\ref{anova.sec1}).

Alternatively, as pointed out by \cite{Zhang_Liang:2013}, one can
also globalize the pointwise $F$-test using the supremum of the
pointwise $F$ test statistic over $\calT$: 
\begin{equation}\label{Fmax.sec1}
F_{\max}=\sup_{t\in \calT} F_n(t). 
\end{equation} 
This $F_{\max}$-test is
somewhat similar to the one suggested in
\cite{Ramsay_Silverman:2005} (p.234) where they used the
squared-root of $F_n(t)$ as their test statistic and used a
permutation-based critical value. In \cite{Zhang_Liang:2013}, the
critical value of the $F_{\max}$-test was obtained via
bootstrapping.  By an application of the $F_{\max}$-test to a real
functional dataset, \cite{Zhang_Liang:2013} realized that the
$F_{\max}$-test may have higher power than their GPF test, particularly 
when the functional data are highly or moderately correlated. In
addition, when applied to an ischemic heart ECG dataset, the
$F_{\max}$-test was significant while the GPF test was not, which
lead to contradictory answers to the relevant question arising from
clinical ischemic heart screening. Therefore, a further study on the
asymptotical and finite sample behaviors, in particular in terms of
level accuracy and power, of the $F_{\max}$-test is necessary.

The major contributions of this paper are as follows. First of all, we derive the asymptotic random expression of the $F_{\max}$ test statistic under very general conditions and the null hypothesis (\ref{anova.sec1}). This allows us to prescribe a parametric bootstrap (PB) method for approximating the asymptotic null distribution of the $F_{\max}$ test statistic. However, one difficulty of this PB method is that it  relies on the asymptotic random expression, which means that it is applicable only when the sample sizes $n_1, \ldots,n_k$ are all large. To overcome this difficulty, we propose a nonparametric bootstrap (NPB) method to approximate the null distribution of $F_{\max}$. Via some extensive simulation studies, we found that the simulated null probability density function (pdf) of the $F_{\max}$ test statistic  is skewed when the functional data are moderately or highly correlated, and the estimated null pdf of the $F_{\max}$ test statistic using the NPB method approximates it reasonably well. Secondly, we show that the NPB $F_{\max}$-test has the correct asymptotic level, and it is root-$n$ consistent i.e. its power tends to 1 under local alternatives that depart from the null hypothesis at the order $n^{-1/2}$. Thirdly,  via some simulation studies,  we show that in general the $F_{\max}$-test outperforms the GPF test of \cite{Zhang_Liang:2013}  in terms of size-controlling, and the $F_{\max}$-test is substantially more powerful than the GPF test  when the functional data are highly or moderately  correlated while the former slightly performs worse than the latter otherwise.  This latter result is intuitively reasonable: when the functional data are highly correlated the GPF test tends to average down the information provided by the data and hence has lower power than the $F_{\max}$-test, whereas when the functional data are less correlated the GPF test tends to summarize more uncorrelated information than the $F_{\max}$-test since the latter takes into account only the information at the maximum of the process $F_n(t)$. Since functional data are usually highly or moderately correlated, the $F_{\max}$-test is therefore preferred to the GPF test in general. This effect also explains the aforementioned contradictory results given by the $F_{\max}$ and GPF tests in the ischemic heart example: the ECG signals are highly correlated and only the $F_{\max}$-test (not the GPF test) detects the significant difference between the ECG signals of the normal and those of stable ischemic heart groups.  Thus the $F_{\max}$-test  result does suggest that it may be feasible to use solely ECG signals in ischemic heart screening, without the aid of more risky stress tests which is required in the current clinical practice. Lastly, we mention that it is straightforward to extend the $F_{\max}$-test to other linear regression models for functional data, including higher-way ANOVA for functional data and functional linear models with functional responses.

The rest of this paper is organized as follows. The main theoretical results on the level accuracy and local power of the $F_{\max}$-test are presented in Section \ref{main.sec2}. In Section \ref{main.sec.discretization} we study the discretization effect on the $F_{\max}$-test, and show that under mild conditions the discretization effect is negligible in terms of both size-controlling and local power. Results of extensive simulation studies and the real data example on ischemic heart screening based on ECG signals are given in Sections \ref{numer.sec} and \ref{section:ischemia:screening} respectively. Some concluding remarks are given in Section \ref{section:conclusion}. Proofs of the theorectical results are deferred to the Appendix.

\section{Main Results}\label{main.sec2}

\subsection{Asymptotic Random Expression under the Null Hypothesis}

We first derive the asymptotic random expression of the $F_{\max}$-test under the null hypothesis (\ref{anova.sec1}). This study  will be helpful for approximating the null distribution of the $F_{\max}$-test.  Notice that for any $t\in\calT$,  the pointwise  between-subject variation $\SSR_n(t)$ given in (\ref{SSRSSE.sec1}) can be expressed as
\begin{equation}\label{SSR.equ}
\SSR_n(t)=\left[\bz_n(t)+\bmu_n(t)\right]^T(\bI_k-\bb_n\bb_n^T/n)\left[\bz_n(t)+\bmu_n(t)\right],
\end{equation} 
 where  $\bI_k$ is the $k\times k$ identity matrix,
 \begin{equation}\label{Zn.sec2}
 \begin{array}{rcl}
 \bz_n(t)&=&\left[\sqrt{n_1}[{\bar{y}}_{1.}(t)-\mu_1(t)],\sqrt{n_2}[{\bar{y}}_{2.}(t)-\mu_2(t)],\cdots,
 \sqrt{n_k}[{\bar{y}}_{k.}(t)-\mu_k(t)]\right]^T,\\
\bmu_n(t)&=&[\sqrt{n_1}\mu_1(t),\sqrt{n_2}\mu_2(t),\cdots,
\sqrt{n_k}\mu_k(t)]^T, \; \\
 \bb_n&=&[\sqrt{n_1},\sqrt{n_2},\cdots, \sqrt{n_k}]^T.
\end{array}
\end{equation} 
Since $\bb_n^T\bb_n/n=1$, it is easy to verify that $\bI_k-\bb_n\bb_n^T/n$ is an idempotent matrix with rank $k-1$. In addition, as $n\rightarrow \infty$, we have
\begin{equation}\label{blim.sec2}
  \bI_k-\bb_n\bb_n^T/n\rightarrow \bI_k-\bb\bb^T,  \;\mbox{ with }
  \bb=[\sqrt{\tau_1},\sqrt{\tau_2},\cdots,
\sqrt{\tau_k}]^T,
\end{equation} 
where $\tau_i=\lim_{n\rightarrow \infty} n_i/n,\; i=1,2,\cdots, k$, are as given in Condition A3 below. Note that  $\bI_k-\bb\bb^T$ in (\ref{blim.sec2}) is also an idempotent matrix of rank $k-1$,  which has the following singular value decomposition:
\begin{equation}\label{bsvd.sec2}
\bI_k-\bb\bb^T=\bU\left(\begin{array}{cc} \bI_{k-1} & \bzero \\
                                          \bzero^T  & 0 \end{array}   \right)\bU^T,
\end{equation} 
where the columns of $\bU$ are the eigenvectors of $\bI_k-\bb\bb^T$.  This property will be used in the derivation of the main results of this paper. Based on the $k$ samples (\ref{ksamp.sec1}),  the pooled
sample covariance function is given by
\begin{equation}\label{hgam.sec2}
 \hgamma(s,t)=(n-k)^{-1} \sum_{i=1}^k \sum_{j=1}^{n_i}[y_{ij}(s)-{\bar{y}}_{i.}(s)][y_{ij}(t)-{\bar{y}}_{i.}(t)],
\end{equation} 
where ${\bar{y}}_{i.}(t), i=1,2,\cdots, k$ are  the sample group mean functions given in  (\ref{ybar.sec1}).

Let $\calL^2(\calT)$ denote the set of all square-integrable functions over $\calT$ and let $C^{\beta}(\calT)$, where $0<\beta\leq1$, denote the set of functions $f$ over $\calT$ which are H\"older continuous with exponent $\beta$ and H\"older modulus $\|f\|_{C^\beta}$. Let $\tr(\gamma)=\int_{\calT} \gamma(t,t) dt$ denote the trace of the covariance function $\gamma(s,t)$.  For the theorectical study, we list  the following regularity conditions.\\

\centerline{\bf Condition A}
\begin{enumerate}
\item The $k$ population mean functions $\mu_1(t), \mu_2(t),\cdots, \mu_k(t)$ in (\ref{ksamp.sec1}) all belong to $\calL^2(\calT)$.

\item The subject-effect functions $v_{ij}(t)=y_{ij}(t)-\mu_i(t), j=1,2,\cdots, n_i; i=1,2, \cdots, k$,  are i.i.d. $\SP(0, \gamma)$ where $0$ denotes the zero function whenever there is no ambiguity.

\item As $n\rightarrow \infty$, the $k$ sample sizes $n_1,\ldots,n_k$
satisfy  $n_i/n \rightarrow \tau_i,\; i=1,2,\cdots, k$, such that
$\tau_1, \tau_2,\cdots, \tau_k \in (0,1)$.

\item  The subject-effect function $v_{11}(t)$ satisfies
$ \E\|v_{11}\|^4=\E\left[\int_{\calT} v_{11}^2(t) dt\right]^2<
\infty. $

\item The covariance function $\gamma(s,t)$ satisfies $\gamma\in C^{2\beta}(\calT\times\calT)$ with $\tr(\gamma)<\infty$, where $0<2\beta\leq 1$. For any $t\in\calT$, $\gamma(t,t)>0$.

\item The expectation $\E\left[ v_{11}^2(s)v_{11}^2(t)\right]$ is uniformly bounded.
That is,   for any $(s,t)\in \calT\times\calT$, we have $\E \left
[v_{11}^2(s)v_{11}^2(t)\right]<C<\infty$ where $C$ is some constant
independent of $(s,t)$.
\end{enumerate}

The first two assumptions A1 and A2 are regular. Note $\calL^2(\calT)$ functions may not be smooth.  Condition A3 requires that the $k$ sample sizes $n_1, n_2, \cdots, n_k$ tend to $\infty$ proportionally. This guarantees that as the total sample size $n\rightarrow \infty$, the sample group mean functions ${\bar{y}}_{i.}(t), i=1,2,\cdots, k$, will converge to some Gaussian processes weakly. The last three assumptions A4, A5 and A6 are imposed so that the pointwise $F$-statistic $F_n(t)$ is well defined at any $t\in\calT$ and the pooled sample covariance function $\hgamma(s,t)$ given in (\ref{hgam.sec2}) converges to the population covariance function $\gamma(s,t)$ uniformly over $\calT\times\calT$.

Let $\GP_k(\bmu,\bGamma)$ denote ``a $k$-dimensional Gaussian process with vector of mean functions $\bmu(t), t\in\calT$, and matrix of covariance functions $\bGamma(s,t), s,t\in\calT$''.  We write $\bGamma(s,t)=\gamma(s,t)\bI_k$ for simplicity when $\bGamma(s,t)=\diag\left[\gamma(s,t), \gamma(s,t),\cdots, \gamma(s,t)\right]$. In particular, $\GP(\eta, \gamma)$ denotes ``a Gaussian process with mean function $\eta(t)$ and covariance function $\gamma(s,t)$''. Further, let ``$\cvgd$'' denote  ``converge in distribution'' in the sense of \cite[p. 474]{Laha_Rohatgi:1979} and \cite[p. 50-51]{vanderVaart_Wellner:1996}. Further,  let ``$X\dequ Y$'' denote that  ``$X$ and $Y$ have the same distribution.'' By Lemma \ref{lem.sec2} in the Appendix, we can  derive the  asymptotical random expression of $F_{\max}$ under the null hypothesis (\ref{anova.sec1}) as given in the following proposition.

\begin{proposition}\label{pro1.sec2}
Under Condition A and the null hypothesis (\ref{anova.sec1}),  as $n\rightarrow \infty$, we have $F_{\max}\cvgd R_0$ with
 \begin{equation}\label{rexp.sec2}
R_0\dequ \sup_{t\in\calT}\Big\{ (k-1)^{-1} \sum_{i=1}^{k-1} w_i^2(t)\Big\},
\end{equation} 
where $w_1,w_2,\cdots, w_{k-1}\iidsim \GP(0,\gamma_w)$ with $\gamma_w(s,t)=\gamma(s,t)/\sqrt{\gamma(s,s)\gamma(t,t)}$.
\end{proposition}

Proposition \ref{pro1.sec2} is useful in discussing a PB method and in investigating an NPB method, detailed in next subsection,  for approximating the null distribution of the $F_{\max}$ test statistic.

\subsection{Approximating the Null Distribution}\label{nullapprox.sec}

By Proposition \ref{pro1.sec2}, $w_1,\cdots, w_{k-1} \iidsim \GP(0,\gamma_w)$ which is known except for $\gamma_w(s,t)$. The covariance function $\gamma_w(s,t)$ can be estimated as
\begin{equation}\label{hgamw.sec2}
\hgamma_w(s,t)=\frac{\hgamma(s,t)}{\sqrt{\hgamma(s,s)\hgamma(t,t)}},
\end{equation} 
where $\hgamma(s,t)$ is  the pooled sample covariance function given in  (\ref{hgam.sec2}). In this case, we can adopt the PB method in \cite{Cuevas_Febrero_Fraiman:2004} to construct an $F_{\max}$-test. The key idea is, for any given  significance level $\alpha$, to  approximate the critical value of the $F_{\max}$-test statistic using the associated upper $100\alpha$-percentile of the random variable $R_0$. Here,  based on (\ref{rexp.sec2}), 
we can sample the i.i.d.  Gaussian processes $w_i, i=1,2,\cdots, k-1$, from $\GP(0, \hgamma_w)$ a large number of times so that a large sample of $R_0$  can be obtained and hence the desired upper $100\alpha$-percentile of $R_0$ can be computed accordingly.  It is obvious that this PB method approximates the critical values well only for large sample sizes.

To overcome the above difficulty of the PB method, we propose an NPB method for
approximating the null distribution of the $F_{\max}$ test statistic
which is applicable for both large and finite sample sizes. Denote
the estimated subject-effect functions as \begin{equation}\label{hvij.sec1}
\hv_{ij}(t)=y_{ij}(t)-{\bar{y}}_{i\cdot}(t),\; j=1,2,\cdots, n_i;
i=1,2,\cdots,k, \end{equation} which can be regarded as estimators of the
subject-effect functions $v_{ij}(t), j=1,2,\cdots, n_{ij},
i=1,2,\cdots, k$, defined in Condition A2. Let
\begin{equation}\label{bksamp.sec2} v_{ij}^{*}(t), j=1,2,\cdots, n_i;
i=1,2,\cdots, k, \end{equation} be bootstrapped $k$ samples randomly generated
from the estimated subject-effect functions given in
(\ref{hvij.sec1}). The nonparametric bootstrapped  $F_{\max}$ test
statistic can then be obtained as \begin{equation}\label{bFmax.sec2}
F^*_{\max}=\sup_{t\in\calT} F_n^*(t), \; \mbox{ where }
F_n^*(t)=\frac{\SSR_n^*(t)/(k-1)}{\SSE^*_n(t)/(n-k)}, \end{equation} with
$\SSR_n^*(t)$ and $\SSE_n^*(t)$ obtained from (\ref{SSRSSE.sec1})
but based on the bootstrapped $k$ samples (\ref{bksamp.sec2}).
Repeat the above bootstrapping process a large number of times and
calculate the upper $100\alpha$-percentile of the bootstrap sample on $F_{\max}^*$, and then
conduct the level-$\alpha$ $F_{\max}$-test using this as the
critical value.

Let $C_{\alpha}$ and $C_{\alpha}^*$ denote the upper $100\alpha$-percentiles of $R_0$ and $F_{\max}^*$ respectively, where $R_0$ is the limit random variable of $F_{\max}$ under the null hypothesis $H_0$ as defined in Proposition \ref{pro1.sec2}.  The following proposition shows that the bootstrapped $F_{\max}$-test statistic $F_{\max}^*$  admits the same limit random expression $R_0$ 
and hence $C_{\alpha}^*$ will also tend to $C_{\alpha}$ in distribution as $n\rightarrow \infty$. Thus the $F_{\max}$-test based on the critical value $C_{\alpha}^*$ has the correct level $\alpha$ asymptotically.

\begin{proposition}\label{pro4.sec2}
Under Condition A,  as $n\rightarrow \infty$, we have $F_{\max}^*\cvgd R_0$ and $C_{\alpha}^*\cvgd C_{\alpha}$.
\end{proposition}

Notice that Proposition \ref{pro4.sec2} holds under both the null and the alternative hypotheses. This is a desirable property since in practice  either of the null and alternative hypotheses may be true.  This property is mainly due to the fact that given the original $k$ samples  (\ref{ksamp.sec1}), the $k$ bootstrapped samples (\ref{bksamp.sec2}) all have the same group mean function as $0$ and hence the null hypothesis always holds for the bootstrapped $k$ samples (\ref{bksamp.sec2}).

Proposition \ref{pro4.sec2} shows that for large samples, the NPB method and the PB method will yield similar approximate critical values for the $F_{\max}$ test statistic. However, the PB method may involve more computational efforts since it requires sampling the Gaussian processes $w_i, i=1,2,\cdots, k-1$, from $\GP(0,\hgamma_w)$ repeatedly. This may not be an easy task when calculation of $\hgamma(s,t)$ is challenging. Besides, larger samples sizes, which may not be available in some applications, are necessary for the PB method to work, in particular when the distribution of $R_0$ is skewed. In general, we prefer the NPB method to the PB method since it can be used under more general conditions than those required by the PB method and its implementation is relatively simple and efficient.

\subsection{The Asymptotic Power}

To  study the asymptotic power of the proposed $F_{\max}$-test, we specify the following local alternative:
\begin{equation}\label{H1n.sec2}
H_{1n} :\;\mu_i(t)=\mu_0(t)+n_i^{-1/2} d_i(t), \; i=1,2,\cdots, k,
\end{equation} 
where $\mu_0(t)$ is the grand mean function, and $d_{1}(t),\cdots,d_{k}(t)$ are any fixed real  functions, independent of $n$. By (\ref{H1n.sec2}), we have $\bmu(t)=[\mu_1(t),\cdots,\mu_k(t)]^T=\mu_0(t)\mathbf{1}_k+[n_1^{-1/2}d_1(t), \cdots, n_k^{-1/2} d_k(t)]^T$ where $\mathbf{1}_k$ denotes the $k\times 1$ vector of ones. It follows that $\bmu_n(t)=\mu_0(t)\bb_n+\bd(t)$ where $\bmu_n(t)$ is defined in (\ref{Zn.sec2}) and $\bd(t)=[d_{1}(t),\cdots,d_{k}(t)]^T$. Under Condition A3, as $n$ tends to $\infty$, the local alternative (\ref{H1n.sec2}) will tend to the null  with the root-$n$ rate. In this sense, if the $F_{\max}$-test can detect the local alternative (\ref{H1n.sec2}) with probability $1$ as long as the information provided by $\bd(t)$ diverges to $\infty$, the $F_{\max}$- test is said to be root-$n$ consistent.  A good test usually should admit the root-$n$ consistency.  Notice that the local alternative (\ref{H1n.sec2}) is the same as the one defined in \cite{Zhang_Liang:2013} who showed that their GPF test is root-$n$ consistent.  In this subsection,  we shall show that the $F_{\max}$-test is also  root-$n$ consistent.

First, since $(\bI_k-\bb_n\bb_n^T/n)\bb_n=\bzero$, under the local alternative (\ref{H1n.sec2}), we have
\begin{equation}\label{SSRexp2.sec2}
\SSR_n(t)=\Big[\bz_n(t)+\bd(t)\Big]^T(\bI_k-\bb_n\bb_n^T/n)\Big[\bz_n(t)+\bd(t)\Big],
\end{equation} 
where $\bz_n(t)$ is defined in (\ref{Zn.sec2}). We then  have the following proposition about the asymptotic distribution of $F_{\max}$ under the local alternative (\ref{H1n.sec2}).

\begin{proposition}\label{pro2.sec2}
Under Condition A and the local alternative (\ref{H1n.sec2}),  as $n\rightarrow \infty$, we have $F_{\max}\cvgd R_1$ with
\begin{equation}\label{rexp1.sec2}
R_1\dequ  \sup_{t\in\calT} \Big\{(k-1)^{-1}\sum_{i=1}^{k-1} [w_i(t)+\delta_i(t)]^2\Big\},
\end{equation} 
where $w_1,\cdots, w_{k-1}\iidsim \GP(0,\gamma_w)$ as in Proposition \ref{pro1.sec2} and $\delta_i(t), i=1,2,\cdots, k-1$, are  the $k-1$ components of $\bdelta(t)=\left(\bI_{k-1}, \bzero\right)\bU^T \bd(t)/\sqrt{\gamma(t,t)}$ with $\bU$ given in (\ref{bsvd.sec2}).
\end{proposition}

With some abuse of  notation, we set $\delta^2=\sum_{i=1}^{k-1} \int_{\calT}\delta_i^2(t) dt$, representing a summary of  the information from $\bd(t)$ with respect to the one-way ANOVA problem (\ref{anova.sec1}) where $\delta_i(t), i=1,2,\cdots, k-1$, are as defined in Proposition ~\ref{pro2.sec2}. By some simple calculations, we can write equivalently $\delta^2=\int_{\calT} \bh(t)^T(\bI_k-\bb\bb^T)\bh(t) dt$ where $\bh(t)=\bd(t)/\sqrt{\gamma(t,t)}$. That is, this $\delta^2$ is the same as the one defined and used in Propositions 2 and 3 in \cite{Zhang_Liang:2013}. By using Proposition 3 in \cite{Zhang_Liang:2013}, we can show the following proposition about the asymptotic power of the NPB $F_{\max}$-test.

\begin{proposition}\label{pro3.sec2}
Under Condition A and the local alternative (\ref{H1n.sec2}), as $n\rightarrow\infty$, the  power of the NPB $F_{\max}$-test $P(F_{\max}\ge C_{\alpha}^*)$  tends to $1$ as $\delta\rightarrow \infty$ where $C_{\alpha}^*$ is the upper $100\alpha$-percentile of the bootstrap test statistic $F_{\max}^*$ defined in (\ref{bFmax.sec2}).
\end{proposition}

Proposition \ref{pro3.sec2}  claims that the  power of the $F_{\max}$-test under the local alternative (\ref{H1n.sec2})  tends to 1 as the information provided by $\bd(t)$ diverges, showing that the proposed $F_{\max}$-test  is root $n$-consistent.  In the proof of Proposition \ref{pro3.sec2}, we shall use the following relationship between the $F_{\max}$-test statistic defined in (\ref{Fmax.sec1}) and the GPF test statistic $T_n$ defined in (\ref{Tn.sec1}):
\[
T_n=\int_{\calT} F_n(t)dt \le (b-a) F_{\max},
\]
where we use the fact that $\calT=[a,b]$.  It then follows that
\begin{equation}\label{Tn2Fmax.sec2}
P(F_{\max}\ge C_{\alpha}^*)\ge P(T_n\ge  (b-a)C_{\alpha}^*).
\end{equation} 
However,  $(b-a)C_{\alpha}^*$ may not be equal or smaller than  the upper $(100\alpha)$-percentile of the GPF test statistic $T_n$. Thus, (\ref{Tn2Fmax.sec2}) does not guarantee that the $F_{\max}$-test has higher power than the GPF test. To compare the powers
of the  two tests, some simulation studies are then needed. Results of such simulation studies are summarized  and discussed  in Section \ref{numer.sec}.

\section{Effect of Discretization on the $F_{\text{max}}$ Test}\label{main.sec.discretization}

We have studied the $F_{\max}$-test in the continuous setup in the previous sections. In practice, accessing the continuous functions $y_{ij}(t)\,, t\in {\calT}$, $j=1,\ldots,n_i$, $i=1,\ldots,k$, may not be always possible, and in most scenarios we have only discretized observations on them. This issue has been partially discussed in \cite{Zhang_Chen:2007} by applying the local polynomial regression (LPR) scheme to reconstruct the continuous functions from the  discretized samples and its asymptotical behavior has been studied therein. Specifically, when the function $\mu_i(t)$ is smooth enough and the discretization points follow some non-degenerate probability distribution function, LPR is applied to reconstruct $\mu_i(t)$. It is shown in \cite{Zhang_Chen:2007} that the LPR estimator is asymptotically unbiased and the convergence rate is discussed. However, in general $\mu_i(t)$ might not be that smooth, and the above approach cannot be applied. In this section we study a more direct approach to approximating the NPB $F_{\max}$-test under discretization and show that the approximation error incurred by discretization is negligible in terms of both the asymptotic level and the asymptotic power under local alternatives.


Suppose we observe the random functions $y_{ij}(t)$, $j=1,\ldots,n_i$, $i=1,\ldots,k$, in (\ref{ksamp.sec1}) only at some discretization points $t_1, \ldots , t_M$ of the interval $\mathcal{T}=[a,b]$, such that $t_1=a$, $t_M=b$ and, for some positive $\tau_M=O(1/M)$ as $M\rightarrow\infty$,
\begin{align}\label{discretization:condition:tau}
0< t_{l+1}-t_{l}\leq \tau_M, \, \mbox{ for all }\, l=1,\ldots,M-1.
\end{align}
Note that the discretization might be non-uniform.
Then we have $k$ samples of random vectors; for $i=1,\ldots,k$, the $i$-th sample consists of
$
 \by_{ij,M}= (y_{ij}(t_1),\ldots,y_{ij}(t_M))^T,\, j=1,\ldots,n_i. \nonumber 
$
 Given the discretized samples $\{\by_{ij,M}\}_{j=1}^{n_i}$, $i=1,\ldots,k$, we have access to the pointwise $F$ test statistics at the discretization points: $F_n(t_l)$, $l=1,\ldots,M$. Thus, parallel to the $F_{\max}$-test 
in the continuous-time case, we can test the discretized null hypothesis
\begin{equation}\label{discrete:H0} H_{0,M}:
\bmu_{1,M}=\ldots=\bmu_{k,M}. 
\end{equation} 
using the following test statistic:
\begin{align}
F_{\max,M}&=\max_{l=1,\ldots,M} F_{n}(t_l).\label{discretization:Fmax_text:definition}
\end{align}
For each $i=1,\ldots,k$, based on the $i$-th discretized sample
$\{\by_{ij,M}\}_{j=1}^{n_i}$, we can generate randomly bootstrapped
sample
$ \{\bv_{ij,M}^{*}\}_{j=1}^{n_i}, $
from the estimated subject-effects vectors$
\{\hbv_{ij,M}=\by_{ij,M}-\bar{\bf y}_{i\cdot,M}\}_{j=1}^{n_i}, $ which
are the estimators of the discretized subject-effects vectors $
\{\bv_{ij,M}= \by_{ij,M} - \bmu_{i,M}\}_{j=1}^{n_i} $ where
$\bar{\bf y}_{i\cdot,M}=n_i^{-1} \sum_{j=1}^{n_i}\by_{ij,M}$ and
$\bmu_{i,M}= (\mu_i(t_1),\ldots,\mu_i(t_M))^T$ is the discretization
of the mean function $\mu_{i}(t)$.
The bootstrapped $F_{\max,M}$ test statistic is then obtained as
\begin{equation}\label{discrete:bFmax.sec2}
F^*_{\max,M}=\max_{l=1,\ldots,M} F_{n}^*(t_l), 
 \end{equation} 
 where $F_{n}^*(t_l)$ are the bootstrap version of $F_{n}(t_l)$ based on the bootstrapped $k$ samples $\{\bv_{ij,M}^{*}\}_{j=1}^{n_i},~ i=1,2,\cdots, k$. 
Repeat the above bootstrapping process a large number of times and
calculate the upper $100\alpha$-percentile of $F_{\max,M}^*$. Then
we can conduct accordingly the NPB $F_{\max,M}$-test for testing the
discretized null hypothesis $H_{0,M}$ (\ref{discrete:H0}).

Let $N_{l}(\boldsymbol{\nu},\boldsymbol{\Gamma})$ denote the distribution of an $l$-vector of Gaussian entries, which has mean vector $\boldsymbol{\nu}$  and covariance matrix  $\boldsymbol{\Gamma}$.  Proposition \ref{Fmax_null:discretization} studies the null distribution of $F_{\max,M}$.

\begin{proposition}
\label{Fmax_null:discretization}
Suppose Condition A holds. Then, under the null hypothesis $H_0$ given in (\ref{anova.sec1}), we have $F_{\max,M}\cvgd R_{0,M}$  as $n\rightarrow \infty$ where
$ \displaystyle{R_{0,M} \dequ \max_{l=1,\ldots,M}}$ $\big\{
(k-1)^{-1} \sum_{i=1}^{k-1} w_{i,M}^2(l)\big\},$
with $\bw_{i,M}\equiv (w_{i,M}(1),\ldots,w_{i,M}(M))^T, i=1,\ldots,k-1$, 
$\iidsim N_{M}(0,\bGamma_{w,M})$ and 
 $\bGamma_{w,M}=\big(\gamma_w(t_p,t_q)\big)_{p,q=1,\ldots,M}$. 
In addition, for any given discretization
$\{t_l\}_{l=1}^M$ of $\calT$ satisfying
(\ref{discretization:condition:tau}),  we may choose $\tbeta$ such
that $0<\tilde{\beta}<\beta$ and  with probability tending to 1, we
have
$
\big|R_0- R_{0,M}\big|\leq \tilde{c} \tau_M^{\tilde{\beta}},
$
where $R_0$ is defined in (\ref{rexp.sec2}) and the constant $\tilde{c}$ depends on the chosen $\tilde{\beta}$.
\end{proposition}

Let $C_{\alpha,M}$ and $C_{\alpha,M}^*$ denote the upper
$100\alpha$-percentiles of $R_{0,M}$ and $F_{\max,M}^*$
respectively. Proposition \ref{discretization:NPB} shows that the
bootstrapped $F_{\max,M}$-test statistic $F_{\max,M}^*$ converges in
distribution to the same $R_{0,M}$ as the $F_{\max,M}$-test
statistic does under $H_0$, and hence $C_{\alpha,M}^*$ will also
tend to $C_{\alpha,M}$ in distribution as $n\rightarrow \infty$.
Thus, the $F_{\max,M}$-test based on the bootstrap critical value
$C_{\alpha,M}^*$ has the correct asymptotic level for testing the
discretized null hypothesis (\ref{discrete:H0}).
In addition, the second result in Proposition \ref{Fmax_null:discretization} 
implies that 
$F_{\max,M}$ and $F_{\max}$ have the same limit test statistic $R_0$
under the null hypothesis, thus the  NPB $F_{\max,M}$-test for
testing  (\ref{discrete:H0}) has the same asymptotic level as the
NPB $F_{\max}$-test for testing $H_0$.

\begin{proposition} \label{discretization:NPB}
Under Condition A, $F_{\max,M}^*\cvgd R_{0,M}$ and $C_{\alpha,M}^*\rightarrow C_{\alpha,M}$ as $n\rightarrow\infty$.
\end{proposition}



Next, we study the local power of the NPB $F_{\max,M}$-test under
the local alternative $H_1$ given in (\ref{H1n.sec2}). The
asymptotic distribution of the $F_{\max,M}$-test statistic under
$H_1$ is given in the following proposition.

\begin{proposition}\label{discretization:pro2.sec2} Suppose Condition A holds. Then, under the local alternative $H_{1}$, 
as $n\rightarrow \infty$, we have $F_{\max,M}\cvgd R_{1,M}$, where
$ R_{1,M}\dequ  \max_{l=1,\ldots,M} \big\{(k-1)^{-1}\sum_{i=1}^{k-1}
[w_{i,M}(l)+\delta_{i,M}(l)]^2\big\}$
$w_{i,M}(l), l=1,\ldots, M$, are defined  in Proposition~\ref{Fmax_null:discretization} 
and $\delta_{i,M}(l), i=1,2,\cdots, k-1$, are  the $k-1$ components of $\bdelta_M(l)=\left(\bI_{k-1}, \bzero\right)\bU^T \bd(t_l)/\sqrt{\gamma(t_l,t_l)}$
with $\bU$ given in (\ref{bsvd.sec2}). Furthermore, suppose $d_i\in C^{\tilde{\beta}}(\calT)$, $i=1,\ldots,k$, for a given $0<\tilde{\beta}<\beta$.
Then, for any given discretization $\{t_l\}_{l=1}^M$ of $\calT$  satisfying (\ref{discretization:condition:tau}), 
with probability tending to 1, 
we have
$
\big|R_1- R_{1,M}\big|\leq \tilde{c} \tau_M^{\tilde{\beta}},
$
where $R_1$ is defined in (\ref{rexp1.sec2}) and the constant
$\tilde{c}$ depends on the chosen $\tilde{\beta}$ and the H\"older
modulus of $d_i(t), i=1,2,\cdots, k$.
\end{proposition}

Combining the results in Propositions
\ref{Fmax_null:discretization}, \ref{discretization:NPB} and
\ref{discretization:pro2.sec2}, as $M\rightarrow\infty$,  we have
that the NPB $F_{\max,M}$-test has the same asymptotic power as the
NPB $F_{\max}$-test under the local alternative $H_1$. Together with
Proposition \ref{pro3.sec2}, this implies that under the local
alternative $H_1$ the  power of the NPB $F_{\max,M}$-test tends to
1, provided that $\bd(t)$ diverges as $n\rightarrow\infty$ and both
Condition A and  (\ref{discretization:condition:tau}) hold.


\section{Simulation Studies}\label{numer.sec}

In this section, we  present some results of simulation studies,  aiming  to
check if the bootstrapped null distribution  of the $F_{\max}$-test
approximates the underlying null distribution of the $F_{\max}$-test well and
how  the $F_{\max}$-test  is compared with  the GPF test of \cite{Zhang_Liang:2013}.

The
$k$ functional samples  (\ref{ksamp.sec1}) were generated from the
following one-way ANOVA model: \begin{equation}\label{simmod.sec3}
\begin{array}{rcl}
y_{ij}(t)&=&\mu_i(t)+v_{ij}(t),\;\;
\mu_i(t)=\bc_i^T[1,t,t^2,t^3]^T,\; \;
 v_{ij}(t)=\bb_{ij}^T\bPsi(t),\; t\in [0,1], \\
\bb_{ij}&=&[b_{ij1},b_{ij2},\cdots,b_{ijq}]^{T}, \; b_{ijr}\dequ \sqrt{\lambda_r}z_{ijr}, \; r=1,2,\cdots,q;
\end{array}
\end{equation} 
$ j=1,2,\cdots, n_i, i=1,2,\cdots, k$, where the parameter vectors
$\bc_i=[c_{i1},c_{i2},c_{i3},c_{i4}]^T, i=1,2,\cdots, k$, for the
group mean functions $\mu_i(t), i=1,2,\cdots, k$, can be flexibly
specified, the random variables $z_{ijr}, r=1,2,\cdots,q; j=1,2,\cdots, n_i; i=1,2,\cdots, k$ are i.i.d. with mean $0$ and variance $1$,  $\bPsi(t)=[\psi_1(t),\cdots, \psi_q(t)]^T$ is a vector of
$q$ orthonormal basis functions  $\psi_r(t), t\in [0,1],
r=1,2,\cdots, q$,  and the variance components $\lambda_r,
r=1,2,\cdots, q$ are positive and decreasing in $r$,
and the number of the basis functions $q$ is an odd positive
integer. These tuning parameters help specify the group mean functions $\mu_i(t)=c_{i1}+c_{i2}t+c_{i3}t^2+c_{i4}t^3, i=1,2,\cdots,k$ and the common covariance function
$
\gamma(s,t)=\bPsi(s)^T\diag(\lambda_1,\lambda_2,\cdots, \lambda_q) \bPsi(t)=\sum_{r=1}^q \lambda_r \psi_r(s)\psi_r(t),
$
 of the subject-effect functions $v_{ij}(t), j=1,2,\cdots, n_i, i=1,2,\cdots, k$.  In the simulations we considered, for simplicity,  we assume that the design time points  for all the functions $y_{ij}(t), j=1,2,\cdots, n_i, i=1,2,\cdots, k$ are the same and are specified as
$t_j=j/(M+1),\ j=1,\cdots, M$, where $M$ is some positive integer. In
practice,  these functions can be observed at different design time
points. In this case,  some smoothing technique, such as the one discussed in \cite{Zhang_Liang:2013}, can be used to reconstruct the functions
$y_{ij}(t), j=1,2,\cdots, n_i,i=1,2,\cdots, k$, and then to evaluate them
at a common grid of time points.   This latter simulation setup will
be time-consuming to carry out and we did not explore it in the simulations.

\begin{figure}
\centerline{\includegraphics[scale=0.70]{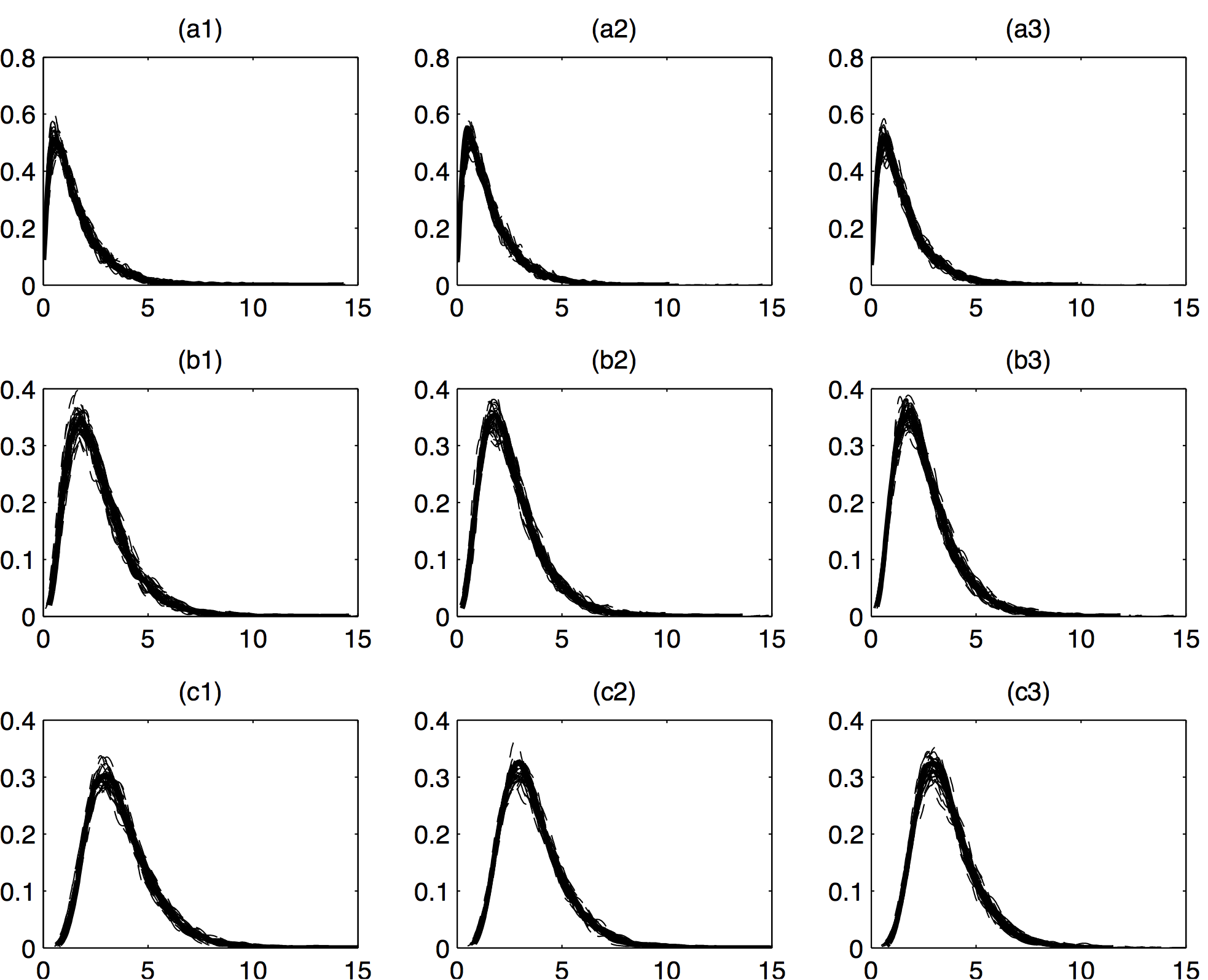}}
\caption{The simulated null  pdfs (wider solid) and the first $50$ bootstrapped null pdfs (dashed) of the $F_{\max}$-test
when  $z_{ijr}, r=1,\cdots, q; j=1,\cdots, n_i; i=1,\cdots, k$, are i.i.d. $N(0,1)$ and $M=80$. From left to right, the columns correspond to the sample size vectors $\bn_1, \bn_2,$ and $\bn_3$ respectively, and from top to bottom the rows correspond to $\rho=0.10, 0.50$ and $0.90$ respectively (those plots associated
 with $\rho=0.30$ and $0.70$ are not shown to save space).} \label{fig1.sim}
\end{figure}

\begin{figure}
\centerline{\includegraphics[scale=0.70]{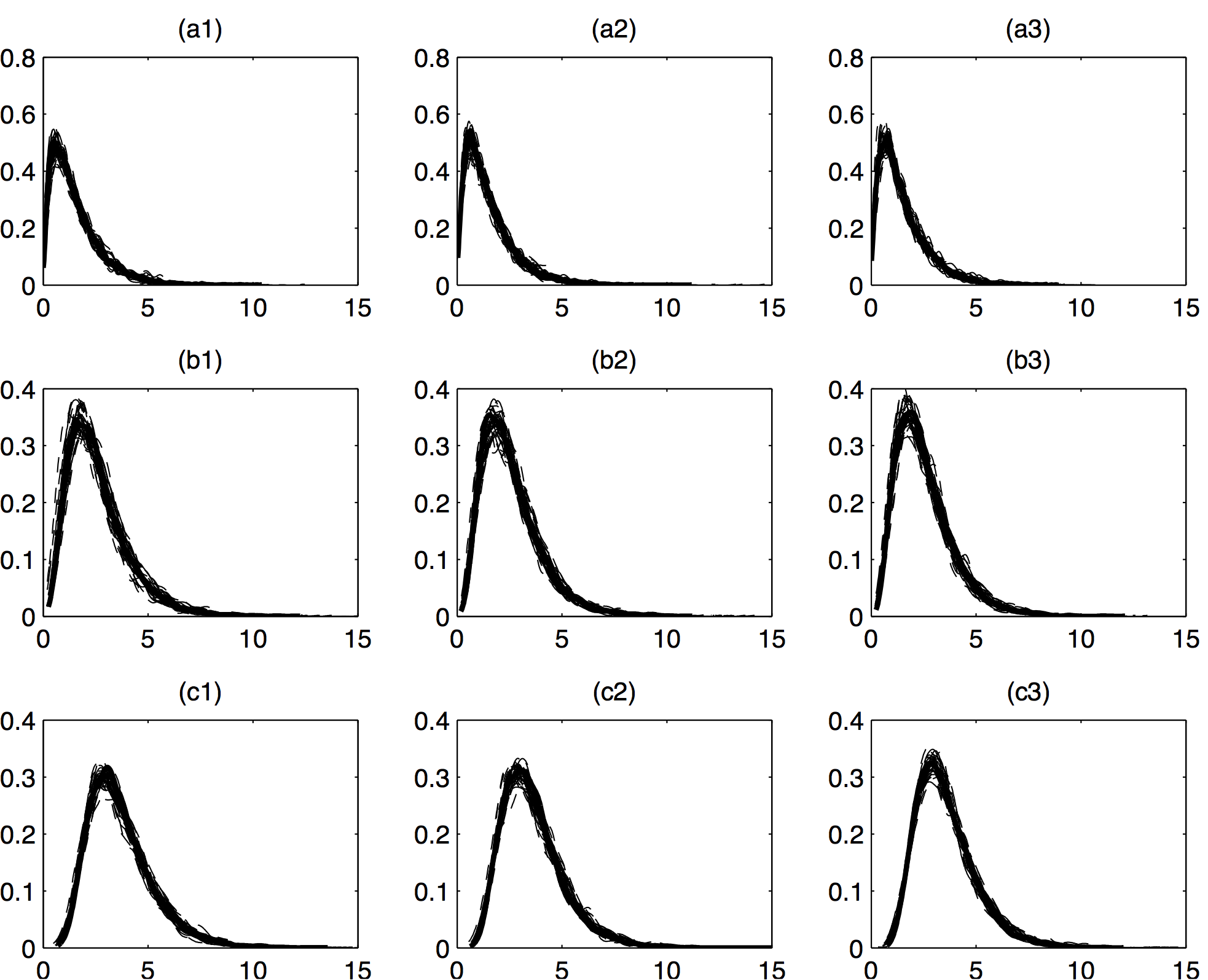}}
\caption{The same as in Fig. \ref{fig1.sim} except that now  $z_{ijr}, r=1,\cdots, q; j=1,\cdots, n_i; i=1,\cdots, k\iidsim t_4/\sqrt{2}$.} \label{fig2.sim}
\end{figure}

We now specify the model parameters in (\ref{simmod.sec3}). To specify the group mean functions $\mu_1(t), \cdots,\mu_k(t)$, 
we set
$\bc_1=[1,2.3,3.4,1.5]^T$ and $\bc_i=\bc_1+(i-1)\delta \bu,
i=2,\cdots, k$, where the tuning parameter $\delta$ specifies the
differences $\mu_i(t)-\mu_1(t), i=2,\cdots, k$, and the constant
vector $\bu$ specifies  the direction of these differences. We set
$\delta$ properly as listed in Tables~\ref{tab1.sim} and ~\ref{tab2.sim} below  so that the null
hypothesis (when $\delta=0$) and the four alternatives (when
$\delta>0$) are fulfilled. In addition, we set
$\bu=[1,2,3,4]^T/\sqrt{30}$ so that it is a unit vector. To specify the common
 covariance function $\gamma(s,t)$, for simplicity, we set $\lambda_r=a \rho^{r-1},
r=1,2,\cdots, q$, for some $a>0$ and $0<\rho<1$.  Notice  that the
tuning parameter $\rho$ not only determines the decay rate of
$\lambda_1,\cdots, \lambda_q$, but also determines how the simulated
functional data are correlated: when $\rho$ is close to $0$,
$\lambda_1,\cdots, \lambda_q$ will decay very fast, indicating that
the simulated functional data are highly correlated; and when $\rho$
is close to $1$, $\lambda_r, r=1,2,\cdots, q$ will decay very slowly,
indicating that the simulated functional data  are nearly uncorrelated.
For simplicity, we set the basis functions as $ \psi_1(t)=1, \;
\psi_{2r}(t)=\sqrt{2}\sin(2\pi r t),\;
\psi_{2r+1}(t)=\sqrt{2}\cos(2\pi r t),\;
  t\in[0,1],$ $ r=1,\cdots, (q-1)/2.$  In addition, we set
$a=1.5, q=11$ and $\rho=0.10, 0.30, 0.50, 0.70, 0.90$ to consider
the five cases when the simulated functional data have  very high,
high,  moderate, low, and very low correlations. We set the number of groups $k=3$ and
specified three cases of the sample size vector $\bn=[n_1,n_2,n_3]$:
$\bn_1=[20,30,30], \bn_2=[40,30,70]$ and $\bn_3=[80,70,100]$,
representing the  small, moderate and larger sample size cases respectively, and
specify two cases of the number of design time points $M$: $M=80$
and $M=150$. Finally, we specified  two cases of the distribution of
the i.i.d. random variables $z_{ijr}, r=1,\cdots, q;\; j=1,\cdots,
n_i; \; i=1,\cdots, k$: $z_{ijr}\iidsim N(0,1)$ and $z_{ijr}\iidsim
t_{4}/\sqrt{2}$, allowing to generate Gaussian and non-Gaussian
functional data respectively with $z_{ijr}$ having mean $0$ and
variance $1$. Notice that the $t_4/\sqrt{2}$ distribution is chosen
since it has nearly the heaviest tails among the $t$-distributions with
finite first two moments.

\begin{table}
\begin{center}
\footnotesize \caption{Empirical sizes and powers (in percentages)
of the GPF and $F_{\max}$  tests for the functional one-way ANOVA
problem (\ref{anova.sec1}) when  the nominal level is $5\%$, $z_{ijr}, r=1,\cdots, q;
j=1,\cdots, n_i; i=1,\cdots, k$, are i.i.d. $N(0,1)$ and $M=80$. The
associated standard deviations (in percentages) are given in
parentheses.}\label{tab1.sim}
\scalebox{0.90}{%
\small\addtolength{\tabcolsep}{-3pt}
\begin{threeparttable}\footnotesize
\begin{tabular}{cc | cc  | cc| cc |cc | cc |cc}\hline\hline
 $\rho$ & $\bn$  & GPF & $F_{\max}$ & GPF  & $F_{\max}$ & GPF
& $F_{\max}$ & GPF & $F_{\max}$ & GPF  & $F_{\max}$ & GPF &$F_{\max}$ \\ \hline

&&\multicolumn{2}{|c}{\underline{$\delta=0$}}&\multicolumn{2}{|c}{\underline{$\delta=0.03$}}
 &\multicolumn{2}{|c}{\underline{$\delta=0.06$}}&\multicolumn{2}{|c}{\underline{$\delta=0.10$}}
 &\multicolumn{2}{|c}{\underline{$\delta=0.13$}}&\multicolumn{2}{|c}{\underline{$\delta=0.16$}}\\ 

 $0.10$&  $\bn_1$ &5.60&  5.04& 7.14& 8.80& 11.20& 26.28&   26.64&   71.96&   43.16&   94.16&   64.62&   99.30\\
&&(0.32)& (0.30)&(0.36)& (0.40)& (0.44)& (0.62)& (0.62)& (0.63)& (0.70)& (0.33)& (0.67)& (0.11)\\
&$\bn_2$&   5.54&    5.09&    8.68&   16.03&   20.82&   58.60&   54.66&   98.06&   81.64&  100.00&    96.02&  100.00\\
&&    (0.32)& (0.31)& (0.39)& (0.51)& (0.57)& (0.69)& (0.70)& (0.19)& (0.54)&         (0.00)   & (0.27)&         (0.00)\\
&$\bn_3$&    5.72&    5.24&   11.20&   26.34&   33.22&   86.64&   83.52&  100.00&    98.44&  100.00&   100.00&   100.00\\
&&    (0.32)& (0.31)& (0.44)& (0.62)& (0.66)& (0.48)& (0.52)&(0.00)    &(0.17)      &(0.00)   &(0.00) &(0.00)\\ \hline

 &&\multicolumn{2}{|c}{\underline{$\delta=0$}}&\multicolumn{2}{|c}{\underline{$\delta=0.05$}}
 &\multicolumn{2}{|c}{\underline{$\delta=0.10$}}&\multicolumn{2}{|c}{\underline{$\delta=0.15$}}
 &\multicolumn{2}{|c}{\underline{$\delta=0.20$}}&\multicolumn{2}{|c}{\underline{$\delta=0.25$}}\\ 

0.30&$\bn_1$& 5.40&5.00&     7.28&    7.24&   11.26&   18.06&   20.88&   40.68&   35.22&   71.54&   52.74&   91.22\\
&&    (0.31)& (0.30)& (0.36)& (0.36)& (0.44)& (0.54)& (0.57)& (0.69)& (0.67)& (0.63)& (0.70)& (0.40)\\
&$\bn_2$&     5.84&    5.02&    8.98&   11.10&   20.22&   39.60&   41.30&   81.49&   70.56&   98.70&   89.94&  100.00\\
&&     (0.33)& (0.30)& (0.40)& (0.44)& (0.56)& (0.69)& (0.69)& (0.54)& (0.64)& (0.16)& (0.42)&         (0.00)\\
&$\bn_3$&    5.36&    4.86&   10.24&   16.58&   31.28&   67.88&   66.72&   97.90&   92.80&  100.00&    99.66&  100.00\\
&&     (0.31)& (0.30)& (0.42)& (0.52)& (0.65)& (0.66)& (0.66)& (0.20)& (0.36)& (0.00)    &(0.08)&         (0.00)\\ \hline

 &&\multicolumn{2}{|c}{\underline{$\delta=0$}}&\multicolumn{2}{|c}{\underline{$\delta=0.10$}}
 &\multicolumn{2}{|c}{\underline{$\delta=0.20$}}&\multicolumn{2}{|c}{\underline{$\delta=0.30$}}
 &\multicolumn{2}{|c}{\underline{$\delta=0.40$}}&\multicolumn{2}{|c}{\underline{$\delta=0.50$}}\\ 
0.50&$\bn_1$& 5.14&    4.40&8.86&    8.52&   20.58&   27.02&   44.26&   65.54&   73.26&   92.64&   90.92&   99.38\\
&&     (0.31)& (0.29)& (0.40)& (0.39)& (0.57)& (0.62)& (0.70)& (0.67)& (0.62)& (0.36)& (0.40)& (0.11)\\
&$\bn_2$&    4.84&    4.40&   13.03&   15.93&   42.82&   64.38&   80.02&   97.24&   98.34&   99.98&   99.96&  100.00\\
&&     (0.30)& (0.29)& (0.47)& (0.51)& (0.69)& (0.67)& (0.56)& (0.23)& (0.18)& (0.01)& (0.02)&         (0.00)\\
&$\bn_3$&    5.54&    4.82&   19.28&   26.82&   67.78&   90.84&   97.82&  100.00&    99.98&  100.00&   100.00&   100.00\\
&&     (0.32)& (0.30)& (0.55)& (0.62)& (0.66)& (0.40)& (0.20)&     (0.00)    &(0.01)&  (0.00) &(0.00) &(0.00)\\ \hline

 &&\multicolumn{2}{|c}{\underline{$\delta=0$}}&\multicolumn{2}{|c}{\underline{$\delta=0.20$}}
 &\multicolumn{2}{|c}{\underline{$\delta=0.40$}}&\multicolumn{2}{|c}{\underline{$\delta=0.60$}}
 &\multicolumn{2}{|c}{\underline{$\delta=0.80$}}&\multicolumn{2}{|c}{\underline{$\delta=1.00$}}\\ 

0.70&$\bn_1$& 5.76&    4.74&   14.34&   12.24&   47.80&   48.84&   86.10&   91.00&    99.10&   99.70&  100.00&    99.98\\
&&     (0.32)& (0.30)& (0.49)& (0.46)& (0.70)& (0.70)& (0.48)& (0.40)& (0.13)& (0.07)&  (0.00)  &  (0.01)\\
&$\bn_2$&    4.86&    4.63&   25.84&   26.16&   83.86&   89.36&   99.82&  100.00&   100.00&   100.00&   100.00&   100.00\\
&&     (0.30)& (0.29)& (0.61)& (0.62)& (0.52)& (0.43)& (0.05)& (0.00) &(0.00)&(0.00)&(0.00) & (0.00)\\
&$\bn_3$&    5.84&    4.94&   42.64&   46.70&   98.54&   99.50&  100.00&   100.00&   100.00&   100.00&   100.00&   100.00\\
&&     (0.33)& (0.30)& (0.69)& (0.70)& (0.16)& (0.09)&(0.00)&(0.00) &(0.00) &(0.00)& (0.00)&(0.00)\\ \hline

 &&\multicolumn{2}{|c}{\underline{$\delta=0$}}&\multicolumn{2}{|c}{\underline{$\delta=0.30$}}
 &\multicolumn{2}{|c}{\underline{$\delta=0.60$}}&\multicolumn{2}{|c}{\underline{$\delta=0.90$}}
 &\multicolumn{2}{|c}{\underline{$\delta=1.20$}}&\multicolumn{2}{|c}{\underline{$\delta=1.50$}}\\ 

0.90&$\bn_1$& 5.14&    5.14&   14.24&   10.82&   55.22&   42.10&   93.50&   85.70&   99.94&   99.36&  100.00&   100.00\\
&&(0.31)& (0.31)& (0.49)& (0.43)& (0.70)& (0.69)& (0.34)& (0.49)& (0.03)& (0.11)    &(0.00)        & (0.00)\\
&$\bn_2$&    5.12&    5.09&   29.82&   23.02&   92.12&   84.62&   99.96&   99.94&  100.00&   100.00&   100.00&   100.00\\
&&     (0.31)& (0.31)& (0.64)& (0.59)& (0.38)& (0.51)& (0.02)& (0.03)& (0.00) &(0.00)&(0.00)&(0.00)\\
&$\bn_3$&    4.90&4.72&   53.44&   41.08&   99.68&   99.30&  100.00&   100.00&   100.00&   100.00&   100.00&   100.00\\
&&     (0.30)& (0.29)& (0.70)& (0.69)& (0.07)& (0.11)&         (0.00)&(0.00)& (0.00)&(0.00)&(0.00)&(0.00)\\ \hline\hline
\end{tabular}
$\bn_1=[20,30,30]$, $\bn_2=[40,30,70]$, $\bn_3=[80,70,100]$.
\end{threeparttable}
}
\end{center}
\end{table}

For a given model configuration, the $k$ functional samples
were generated as in (\ref{simmod.sec3}) . Then the
GPF test statistic and its P-value were computed.
At the same time, the $F_{\max}$ test statistic was computed, and $10000$
bootstrap replicates were generated to compute the P-value of the
NPB $F_{\max}$-test. When the P-value of a testing procedure
is smaller than the nominal significance level $\alpha$ ($5\%$
here), the null hypothesis (\ref{anova.sec1}) is rejected. The above
process was repeated $N=5000$ times. The empirical size or power of
a testing procedure was then computed as the proportion of the number
of rejections to the number of replications $N=5000$.

We first check if the bootstrapped null pdf of the $F_{\max}$ test statistic
works well in approximating the underlying null pdf of the
$F_{\max}$ test statistic. To this end, each panel in Fig.~\ref{fig1.sim}
displays the simulated null  pdf (wider solid curve) and the first
$50$ bootstrapped null pdfs (dashed curves) of the $F_{\max}$ test statistic
when $z_{ijr}, r=1,\cdots, q; j=1,\cdots, n_i; i=1,\cdots, k\iidsim
N(0,1)$ and $M=80$. 
Each of the pdfs was computed using the usual kernel density estimator \citep[Ch. 2]{Wand_Jones:1995} based on the $5000$ simulated $F_{\max}$ test statistics
(\ref{Fmax.sec1}) or the $10000$ bootstrapped
$F_{\max}$ test statistic  when $\delta=0$ and $M=80$. From
Fig.~\ref{fig1.sim}, it is seen that the bootstrapped null pdfs of
the $F_{\max}$ test statistic approximate the associated simulated null pdf
rather well in all of the nine panels, showing that the NPB method does
work reasonably well in approximating the underlying null pdfs of
the $F_{\max}$ test statistic.  Furthermore, it seems  that the effects of the
sample sizes are  not remarkable; but the shapes of
 the simulated and bootstrapped null pdfs of the $F_{\max}$ test statistic are
strongly affected  by the decay rates  of the  variance components
$\lambda_r, r=1,2,\cdots,q$, namely, stronger correlation in the functional data ($\rho$ smaller) causes
more skewness in the null distribution of $F_{\max}$.

In Fig.~\ref{fig2.sim}, we display the simulated and bootstrapped
null pdfs of the  $F_{\max}$-test when $z_{ijr}, r=1,\cdots, q;
j=1,\cdots, n_i; i=1,\cdots, k\iidsim t_4/\sqrt{2}$. From
Fig.~\ref{fig2.sim}, we can  see that  the NPB method still works
reasonably well in approximating the underlying null pdf of the
$F_{\max}$-test and  that the shapes of
 the simulated and bootstrapped null pdfs of the $F_{\max}$-test  are
again strongly affected  by the decay rates  of the  variance components
$\lambda_r, r=1,2,\cdots,q$. We  also observe that the effects of the
sample sizes  are comparably minor.

\begin{table}
\begin{center}
\footnotesize \caption{Empirical sizes and powers (in percentage) of
the GPF and $F_{\max}$-tests for the functional
one-way ANOVA problem (\ref{anova.sec1}) when the nominal size is $5\%$,  $z_{ijr}, r=1,\cdots,
q; j=1,\cdots, n_i; i=1,\cdots, k$, are i.i.d. $t_4/\sqrt{2}$ and $M=80$. The associated
standard deviations (in percentage) are given in
parentheses.}\label{tab2.sim}
\scalebox{0.90}{%
\small\addtolength{\tabcolsep}{-3pt}
\begin{threeparttable}\footnotesize
\begin{tabular}{ c c | cc  | cc| cc |cc | cc |cc}\hline\hline
 $\rho$ & $\bn$  & GPF & $F_{\max}$ & GPF  & $F_{\max}$ & GPF
& $F_{\max}$ & GPF & $F_{\max}$ & GPF  & $F_{\max}$ & GPF &$F_{\max}$ \\ \hline
&&\multicolumn{2}{|c}{\underline{$\delta=0$}}&\multicolumn{2}{|c}{\underline{$\delta=0.03$}}
 &\multicolumn{2}{|c}{\underline{$\delta=0.06$}}&\multicolumn{2}{|c}{\underline{$\delta=0.10$}}
 &\multicolumn{2}{|c}{\underline{$\delta=0.13$}}&\multicolumn{2}{|c}{\underline{$\delta=0.16$}}\\ 

0.10&$\bn_1$&      5.42&    5.00&     6.51&    8.93&   12.30&   29.94&   28.46&   74.80&   46.08&   93.24&   65.12&   98.66\\
&&    (0.32)& (0.30)& (0.34)& (0.40)& (0.46)& (0.64)& (0.63)& (0.61)& (0.70)& (0.35)& (0.67)& (0.16)\\
&$\bn_2$&    5.26&    4.80&8.00&    15.24&   21.92&   62.02&   55.92&   97.34&   82.56&   99.86&   95.64&   99.98\\
&&    (0.31)& (0.30)& (0.38)& (0.50)& (0.58)& (0.68)& (0.70)& (0.22)& (0.53)& (0.05)& (0.28)& (0.01)\\
&$\bn_3$&    5.08&    5.00&    11.26&   27.60&   35.82&   86.94&   83.58&   99.92&   98.04&  100.00&    99.70&  100.00\\
&&    (0.31)& (0.30)& (0.44)& (0.63)& (0.67)& (0.47)& (0.52)& (0.03)& (0.19)&(0.00)&(0.07)&(0.00)\\ \hline

    &&\multicolumn{2}{|c}{\underline{$\delta=0$}}&\multicolumn{2}{|c}{\underline{$\delta=0.05$}}
 &\multicolumn{2}{|c}{\underline{$\delta=0.10$}}&\multicolumn{2}{|c}{\underline{$\delta=0.15$}}
 &\multicolumn{2}{|c}{\underline{$\delta=0.20$}}&\multicolumn{2}{|c}{\underline{$\delta=0.25$}}\\ 

0.30&$\bn_1$&    5.22&    4.82&    7.26&    7.42&   11.98&   19.28&   22.28&   44.74&   38.52&   74.06&   55.78&   91.16\\
&&    (0.31)& (0.30)& (0.36)& (0.37)& (0.45)& (0.55)& (0.58)& (0.70)& (0.68)& (0.61)& (0.70)& (0.40)\\
&$\bn_2$&    5.14&    4.72&    7.98&   11.04&   20.98&   42.60&   43.58&   81.66&   70.42&   97.72&   89.78&   99.74\\
&&    (0.31)& (0.29)& (0.38)& (0.44)& (0.57)& (0.69)& (0.70)& (0.54)& (0.64)& (0.21)& (0.42)& (0.07)\\
&$\bn_3$&    5.90&5.74&   10.30&   16.92&   32.56&   70.18&   67.42&   97.64&   93.08&   99.94&   99.02&  100.00\\
&&    (0.33)& (0.32)& (0.42)& (0.53)& (0.66)& (0.64)& (0.66)& (0.21)& (0.35)& (0.03)& (0.13)&(0.00)\\ \hline

    &&\multicolumn{2}{|c}{\underline{$\delta=0$}}&\multicolumn{2}{|c}{\underline{$\delta=0.10$}}
 &\multicolumn{2}{|c}{\underline{$\delta=0.20$}}&\multicolumn{2}{|c}{\underline{$\delta=0.30$}}
 &\multicolumn{2}{|c}{\underline{$\delta=0.40$}}&\multicolumn{2}{|c}{\underline{$\delta=0.50$}}\\ 

0.50&$\bn_1$&    4.72&    5.06&    9.10&9.03&   21.36&   29.62&   48.60&   67.40&   73.34&   92.54&   91.46&   98.96\\
&&    (0.29)& (0.30)& (0.40)& (0.40)& (0.57)& (0.64)& (0.70)& (0.66)& (0.62)& (0.37)& (0.39)& (0.14)\\
&$\bn_2$&    4.97&    4.78&   13.80&   17.08&   44.24&   65.58&   80.64&   96.36&   97.58&   99.84&   99.72&   99.98\\
&&    (0.30)& (0.30)& (0.48)& (0.53)& (0.70)& (0.67)& (0.55)& (0.26)& (0.21)& (0.05)& (0.07)& (0.01)\\
&$\bn_3$&    5.66&    5.00&    19.70&   28.16&   68.40&   91.02&   97.46&   99.94&   99.88&  100.00&   100.00&   100.00\\
&&    (0.32)& (0.30)& (0.56)& (0.63)& (0.65)& (0.40)& (0.22)& (0.03)& (0.04)&(0.00)&(0.00)&(0.00)\\ \hline

    &&\multicolumn{2}{|c}{\underline{$\delta=0$}}&\multicolumn{2}{|c}{\underline{$\delta=0.20$}}
 &\multicolumn{2}{|c}{\underline{$\delta=0.40$}}&\multicolumn{2}{|c}{\underline{$\delta=0.60$}}
 &\multicolumn{2}{|c}{\underline{$\delta=0.80$}}&\multicolumn{2}{|c}{\underline{$\delta=1.00$}}\\ 

0.70&$\bn_1$&    4.76&    4.92&   13.10&   12.74&   48.28&   50.82&   85.38&   90.92&   98.76&   99.66&   99.74&  100.00\\
&&    (0.30)& (0.30)& (0.47)& (0.47)& (0.70)& (0.70)& (0.49)& (0.40)& (0.15)& (0.08)& (0.07)&(0.00)\\
&$\bn_2$&    4.86&    4.48&   27.26&   26.42&   84.76&   89.74&   99.50&   99.92&   99.98&  100.00&    99.98&  100.00\\
&&    (0.30)& (0.29)& (0.62)& (0.62)& (0.50)& (0.42)& (0.09)& (0.03)& (0.01)&(0.00)&(0.01)&(0.00)\\
&$\bn_3$&    5.40&5.00&    44.02&   47.70&   98.46&   99.48&  100.00&   100.00&   100.00&   100.00&   100.00&   100.00\\
 &&   (0.31)& (0.30)& (0.70)& (0.70)& (0.17)& (0.10)&(0.00)&(0.00)&(0.00)&(0.00)&(0.00)&(0.00)\\ \hline

&&\multicolumn{2}{|c}{\underline{$\delta=0$}}&\multicolumn{2}{|c}{\underline{$\delta=0.30$}}
 &\multicolumn{2}{|c}{\underline{$\delta=0.60$}}&\multicolumn{2}{|c}{\underline{$\delta=0.90$}}
 &\multicolumn{2}{|c}{\underline{$\delta=1.20$}}&\multicolumn{2}{|c}{\underline{$\delta=1.50$}}\\ 

0.90&$\bn_1$&    4.50&5.40&   14.36&   11.18&   55.52&   43.68&   93.54&   86.90&   99.50&   99.26&   99.94&   99.98\\
 &&   (0.29)& (0.31)& (0.49)& (0.44)& (0.70)& (0.70)& (0.34)& (0.47)& (0.09)& (0.12)& (0.03)& (0.01)\\
 &$\bn_2$&   4.76&    5.08&   29.92&   23.16&   91.68&   85.52&   99.92&  100.00&   100.00&   100.00&   100.00&   100.00\\
  &&  (0.30)& (0.31)& (0.64)& (0.59)& (0.39)& (0.49)& (0.03)&(0.00)&(0.00)&(0.00)&(0.00)&(0.00)\\
  &$\bn_3$&  4.78&    4.74&   51.90&   41.56&   99.50&   99.32&  100.00&   100.00&   100.00&   100.00&   100.00&   100.00\\
  &&  (0.30)& (0.30)& (0.70)& (0.69)& (0.09)& (0.11)& (0.00)&(0.00)&(0.00)&(0.00)&(0.00)&(0.00)\\ \hline\hline
\end{tabular}
$\bn_1=[20,30,30]$, $\bn_2=[40,30,70]$, $\bn_3=[80,70,100]$.
\end{threeparttable}
}
\end{center}
\end{table}

We now turn to check how  the NPB $F_{\max}$-test is compared with the
GPF test of \cite{Zhang_Liang:2013} for the one-way ANOVA problem
(\ref{anova.sec1})  in terms of level accuracy and power.  For
this purpose, we summarize in Tables~\ref{tab1.sim} and ~\ref{tab2.sim}  the
empirical sizes and powers (in percentages) of the GPF and
$F_{\max}$-tests when  $z_{ijr}, r=1,\cdots, q; j=1,\cdots, n_i;
i=1,\cdots, k$ are i.i.d. $N(0,1)$  or  $t_4/\sqrt{2}$, respectively.
The associated standard deviations (in percentages) of these
empirical sizes and powers
 are also given  in the parentheses. From the columns associated with $\delta=0$ in  both tables,
 we see that in terms of
 size controlling, in general, the NPB $F_{\max}$-test outperforms the GPF test.
 From the columns associated with $\delta>0$ in both tables,  we also see that
 in terms of power,  in general the NPB $F_{\max}$-test has higher
 powers than the GPF test except when the functional data are less correlated
 ($\rho=0.90$) and the advantages of the NPB $F_{\max}$-test over the GPF test
 become more significant as the correlation in the functional data  increases.

In the above, we only presented the simulation results when $M=80$
since those when $M=150$ are similar.  From all
these simulation studies,   we  conclude that the NPB method
presented in this paper works reasonably well in approximating the
underlying null pdf of the $F_{\max}$-test and in general the NPB $F_{\max}$ test statistic
 outperforms the GPF test in terms of size-controlling and
it  generally has significantly higher powers than the GPF test when the
functional data are moderately or strongly correlated.

\section{Screening Myocardial Ischemia Using Resting Electrocardiac Signals}\label{section:ischemia:screening}

In this section we present an example where the NPB $F_{\max}$ is useful in discovering new methods in  medical studies. Specifically,
detecting myocardial ischemia (MI) is an important clinical mission, in particular in outpatient evaluations. Typical evaluations include resting surface electrocardiography (ECG) examination and non-invasive stress testing, for example, exercise ECG and the single-photon emission computed tomography (SPECT) thallium scan. Although resting ECG is a fundamental measure in this regard, compared with stress tests, it is not accurate enough in  detecting MI in a typical chest pain clinic \citep{Gibbons_Abrams_Chatterjee:2003}. On the other hand, while stress tests are more accurate than a resting ECG, they have limitations in clinics. Among the limitations, the most important one is the stress itself:  it is associated with the risk of acute attacks during the testing. Thus, stress tests cannot be performed on patients who are extremely susceptible to the provocations used in the tests, and such patients have to directly undergo invasive tests \citep{Gibbons_Abrams_Chatterjee:2003}.

Recent findings from biophysics and pathology interdisciplinary work \citep{Swan:1979}, researches on the physical characteristics of myocardial strain in echocardiography \citep{Urheim_Edvardsen_Torp:2000} and studies on electric signals \citep{delRio_McConnell_Clymer:2005} have necessitated a reappraisal of the ECG information for detecting MI.
For example, spectral analysis of the resting ECG signals from dogs revealed a shift from high- to low-frequency ranges in ischemia cases \citep{Mor-Avi_Akselrod:1990}. Similar phenomena were associated with balloon inflations during percutaneous transcatheter angioplasty in MI patients \citep{Abboud_Cohen_Selwyn:1987,Pettersson_Pahlm_Carro:2000}. Based on these findings, we hypothesize that the ischemic myocardium information is contained in the resting ECG signals with the proper manipulation, and the power spectrum of the manipulated ECG signal is shifted to the low frequency.
To test this hypothesis, the following prospective clinical study was conducted and the proposed NPB $F_{\text{max}}$ test was applied to assess the result.

\subsection{Study populations and categorization}
We designed this study to use the same procedures employed in daily clinical practice and enrolled $393$ consecutive individuals who visited an outpatient clinic complaining of chest pain. We evaluated the possibility of myocardial ischemia based on published guidelines \citep{Gibbons_Balady_Beasley:2002,Gibbons_Abrams_Chatterjee:2003,Klocke_Baird_Lorell:2003}. The institutional review board of the hospital approved this protocol. Patients completed all the examinations within a span of three months, and all provided written consent.

\begin{figure}
\includegraphics[width=0.6\textwidth]{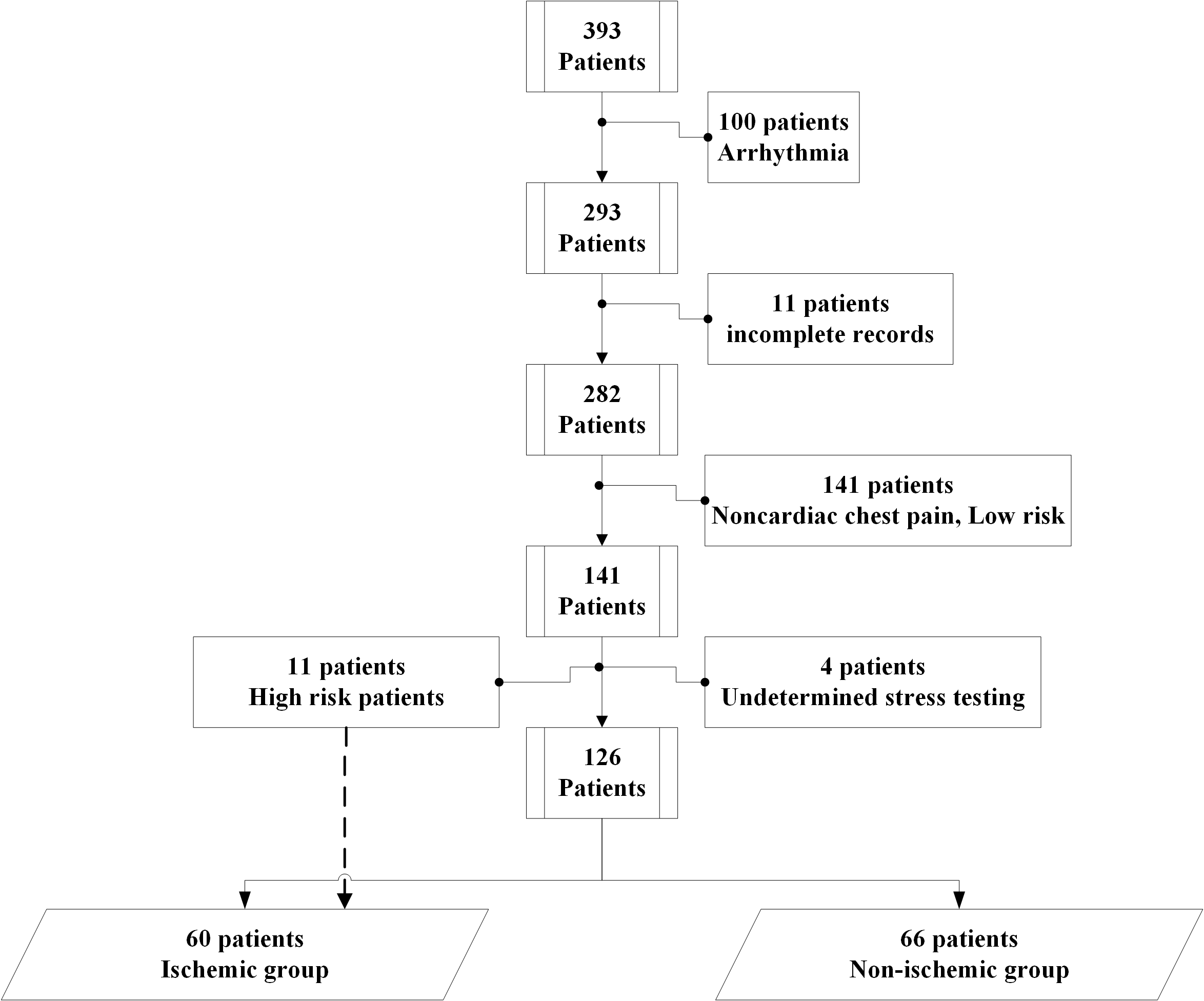}
\caption{Flow chart of patient categories. Note that all $11$ high-risk patients had had a positive CAG and remained in the ischemic group. In the end, there are $71$ patients in the ischemic group and $66$ patients in the normal group.}
\label{fig2}
\end{figure}

We studied the following two groups -- an ischemic group and a normal group defined by the following procedure by taking the clinical procedures into account. The ischemic group consisted of patients who directly underwent a coronary arteriography (CAG) and diagnosed as MI
and patients having positive exercise ECG and/or SPECT thallium scan. If a patient's exercise ECG was inconclusive, the result of SPECT thallium scan was used for classification. The normal group comprised of patients who had had at least one stress test, and any of the test results was negative.
We excluded individuals based on any of the following criteria: (i) non-cardiac chest pain and a very low risk of developing MI; (ii) non-NSR (normal sinus rhythm); (iii) already underwent exercise ECG, but were non-diagnostic or intolerant of the SPECT thallium scan, and unwilling to undertake further evaluation; (iv) a heart rate of less than $40$ beats per minute.
In the end, $137$ patients were eligible for the analysis -- $71$ patients in the ischemic group and 66 patients in the normal group. Please see Fig. \ref{fig2} for the details. 

\subsection{Data acquisition}

We acquired from each subject $85$ seconds of $12$-lead resting ECG signals of his/her first visit to the outpatient clinic. The signals were acquired at $500$Hz and quantized at $12$ bits across $\pm10$ mV \citep{Bailey_Berson_Garson:1990}. The signals were passed through a digital low-pass filter with a $-3$dB cutoff at $60$ Hz, and stored in double precision. Besides removing the inherent amplitude deviations, such as power line noise, thermal noise, etc., we preserved all possible physiological activities.
The recorded $12$ lead ECG signals for the $i$-th subject is denoted as a $12\times {\calN}$ matrix $\boldsymbol{e}^{(i)}$, ${\calN}=85\times 500$, so that $\boldsymbol{e}^{(i)}(l,j)$ is the $l$-th lead ECG signal sampled at time $j\tau$, where $\tau=1/500$ is the sampling interval and $j=1,\ldots,{\calN}$.

\begin{figure}
\includegraphics[width=0.95\textwidth]{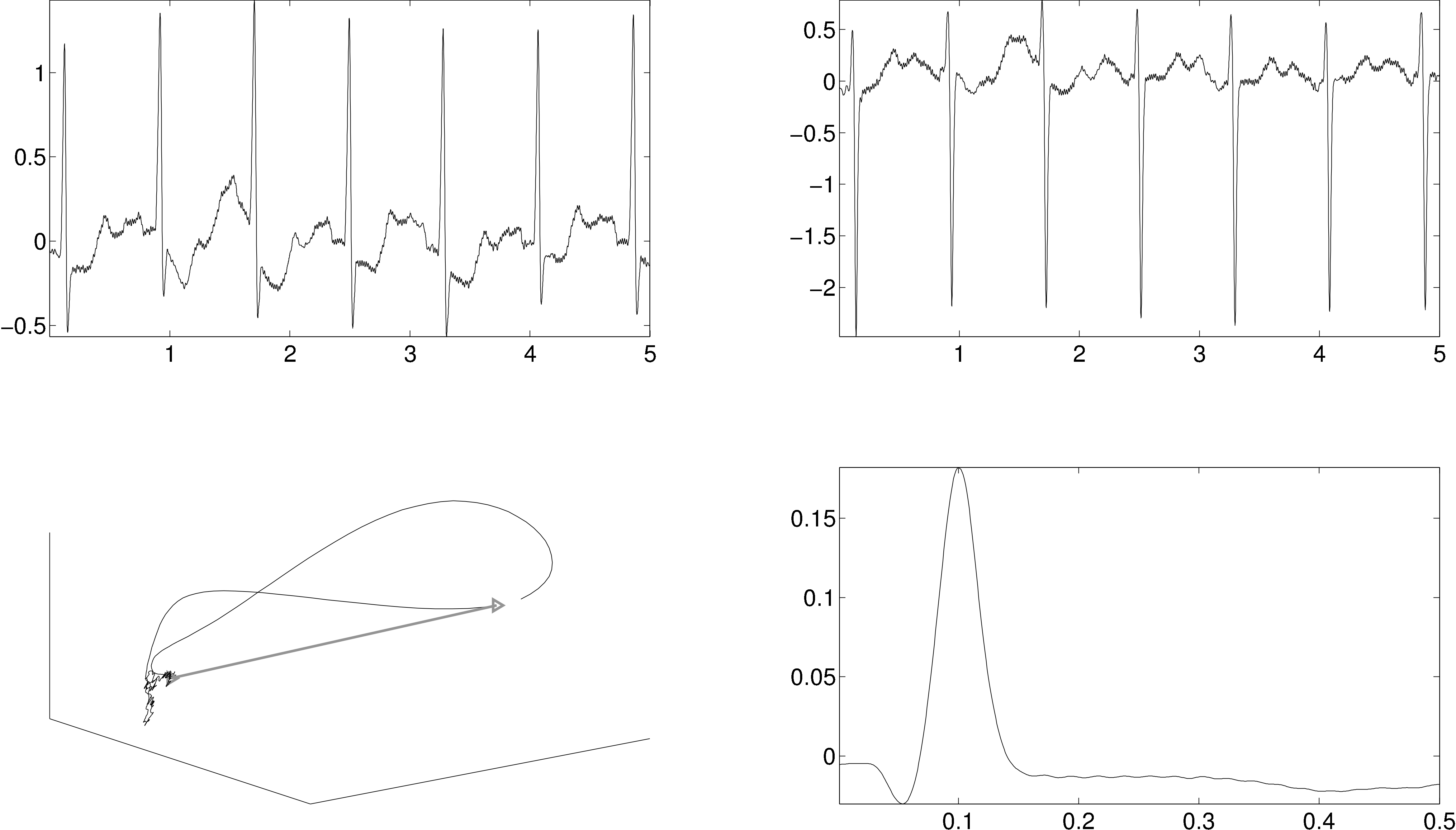}
\caption{Upper left (right): a 5-second
Lead II (Lead III) ECG signal. 
Lower left: the VCG signal of an extracted sinus heartbeat
$\boldsymbol{b}^{(i,l)}$ (black curve) superimposed with the direction
of $\bv^{(i,l)}$ pointing from the origin (i.e. the
atrio-ventricular node) to the R-peak (gray curve). Lower right: the
adaptive lead signal $\bs_0^{(i,l)}$ constructed by
projecting $\boldsymbol{b}^{(i,l)}$ onto $\bv^{(i,l)}$.}
\label{fig1}
\end{figure}

First, for the $i$-th subject, we reduced the influence of the lead system by recovering the 3-dimensional vectocardiograph (VCG) data
with LevkovÕs algorithm \citep{Levkov:1987}, which is denoted as a $3\times {\calN}$ matrix $V^{(i)}$, where the $j$-th column is the dipole current direction in $\RR^3$ at time $j\tau$.

Second, we reduced the heart rate variability (HRV) effect. Denote
$t^{(i,l)}\in\{1,\ldots, {\calN}\}$ to be the indices of time
stamps of the $l$-th R-peak of the $n$-th subject. The R-peaks were
detected using the standard method
\citep{Clifford_Azuaje_McSharry:2006}. We then extract the first $L$
sinus heartbeat signals, where $L\in\mathbb{N}$. Here the $l$-th
sinus heartbeat signal is the VCG signal between time
$t^{(i,l)}\tau$ and $t^{(i,l+1)}\tau$. Denote the $l$-th extracted
sinus heartbeat as a $3\times (t^{(i,l+1)}-t^{(i,l)}+1)$ matrix
$\boldsymbol{b}_0^{(i,l)}$, where
$\boldsymbol{b}_0^{(i,l)}(j,m)=\boldsymbol{e}^{(i)}(j,t^{(i,l)}+m-1)$,
$j=1,\ldots,3$ and $m=1,\ldots,t^{(i,l+1)}-t^{(i,l)}$. Next, we
interpolate each heartbeat by the cubic spline interpolation to be
of the uniform length $500$ to eliminate the HRV influence, and
denote the $l$-th interpolated heartbeat as a $3\times 500$ matrix
$\boldsymbol{b}^{(i,l)}$. Then we ``co-axial project'' each
interpolated heartbeat $\boldsymbol{b}^{(i,l)}$ onto the R-peak
direction, which by our construction is the first column of
$\boldsymbol{b}^{(i,l)}$. Denote this R-peak direction as a column
vector $\bv^{(i,l)}\in\RR^3$, and we project
$\boldsymbol{b}^{(i,l)}$ onto $\bv^{(i,l)}$:
$$
\bs_0^{(i,l)} =[\bv^{(i,l)}]^T\boldsymbol{b}^{(i,l)} \in\RR^{1\times
N}.
$$
We call $\bs_0^{(i,l)}$ an {\it adaptive lead signal}, where the
adjective {\it adaptive} means that the cardiac axis and lead system
effects are reduced.

Third, we eliminated other physiological effects, in particular the
respiration, by normalizing each beat $\bs_0^{(i,l)}$ such that its
$L^2$ norm is $1$, which is denoted as another row vector
$\bs^{(i,l)}$ of length $500$. With $\{\bs^{(i,l)} \}_{l=1}^{L}$,
the {\it adaptive ECG waveform} for the $i$-th subject is defined as
\begin{align}\label{definition:power_spectrum}
\bs^{(i)}=\frac{1}{L}\sum_{l=1}^{L} \bs^{(i,l)}.
\end{align}
We tested our null hypothesis based on the adaptive ECG waveform collected from the subjects. A typical ECG signal, the VCG signal and the adaptive ECG waveform representation are demonstrated in Fig. \ref{fig1}.

\begin{figure}
\begin{center}
\includegraphics[width=1\textwidth]{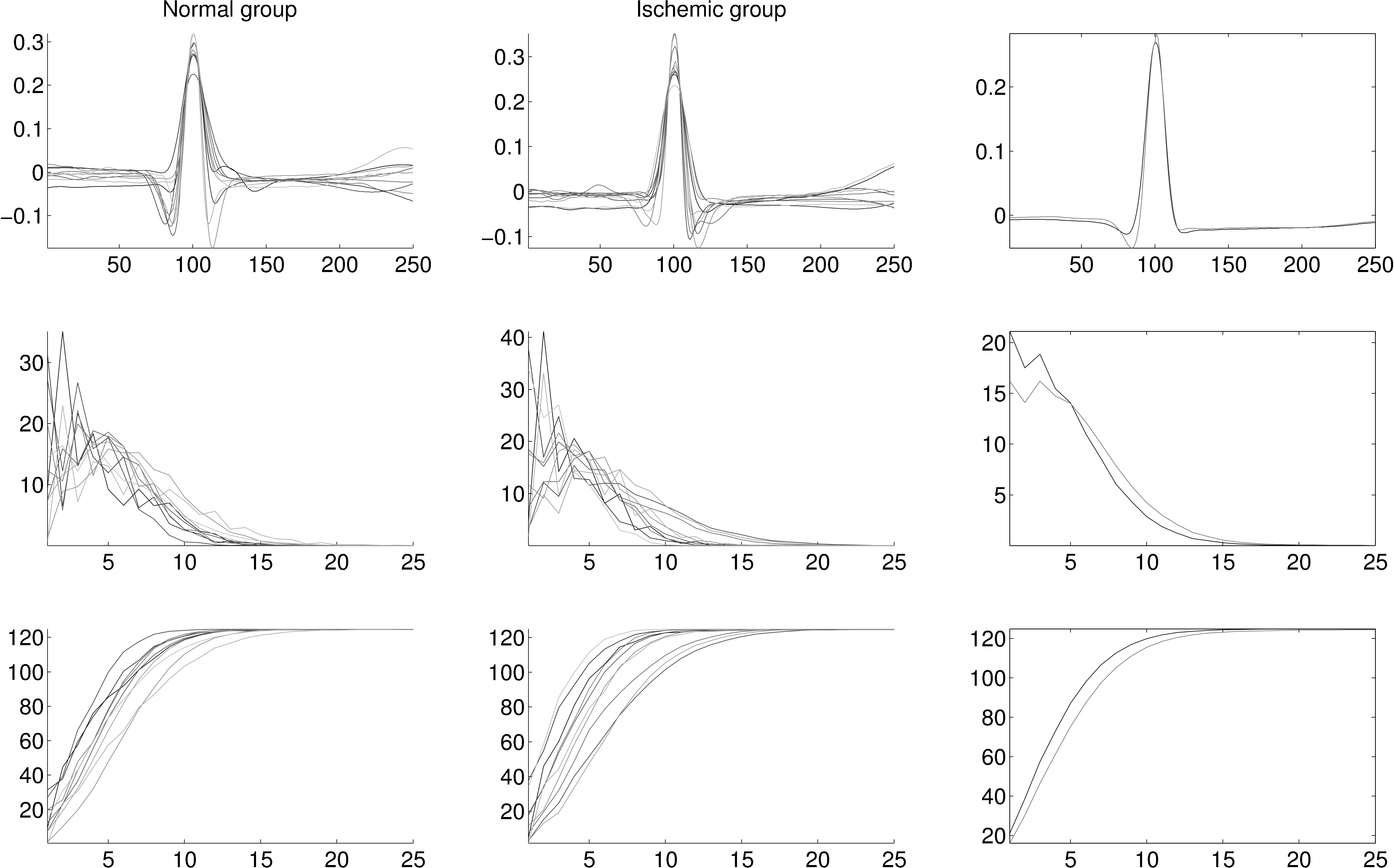}  
\footnotesize\caption{Upper left (middle): ten adaptive ECG waveform representations from ten randomly selected subjects in the normal (ischemia) group. We only show the first $250$ points for the demonstration purpose. 
Upper right: the mean adaptive ECG waveform representations of the normal and ischemia groups are respectively plotted in gray and black.
Note that a notch on the left-hand side of the peak can be visually observed.
Middle left (middle): ten power spectra from ten random subjects in the normal (ischemia) group. We only show the first $25$ Fourier modes for the demonstration purpose. 
Middle right: the mean power spectra of the normal and ischemia groups are respectively plotted in gray and black
Note the obvious shifting of the power spectrum from the high frequency area to the low frequency area in the ischemic group.
Lower left (middle): cumulative power spectra from ten random subjects in the normal (ischemia) group. We only show the first $25$ cumulative powers for the demonstration purpose. 
Lower right: the mean cumulative power spectra of the normal and ischemia groups are respectively plotted in gray and black. }\label{fig3}
\end{center}
\end{figure}

To confirm our hypothesis derived from the results in
\cite{Mor-Avi_Akselrod:1990,Abboud_Cohen_Selwyn:1987,Pettersson_Pahlm_Carro:2000},
that is, the power spectrum of a beat of an ischemic heart is
shifted to the low frequency region, we study the power spectrum of
$\bs^{(i)}$, which is defined as $P^{(i)}=|\mathcal{F} \bs^{(i)}
|^2$, where $\mathcal{F}$ is the discrete Fourier transform. Here
$P^{(i)}$ is a row vector of length $500$. We follow the convention
and say the first entry of $P^{(i)}$ is the DC term, the second to
$250$ entries are of the positive frequencies and the left are of
the negative frequencies. Further, we define the cumulative power
spectrum of $\bs^{(i)}$ as $f^{(i)}(m)=\sum_{j=2}^{m+1} P^{(i)}(j)$,
where $m=1,\ldots, 250$, that is $f^{(i)}$ is a row vector of length
$250$. Fig. \ref{fig3} shows a set of adaptive ECG waveform
$\bs^{(i)}$ (upper row), power spectra $P^{(i)}$ (middle row) and
the electrophysiological fingerprints $f^{(i)}$ (lower row).

\subsection{$F_{\max}$ Analysis and Results}
To test the one-way ANOVA problem (\ref{anova.sec1}) for the stable
cardiac ischemia data, we considered the one-way ANOVA problem over
the adaptive lead signal $\bs_0^{(i,l)}$, the adaptive ECG waveform
$\bs^{(i)}$, its power spectrum $P^{(i)}$ and the cumulative power
spectrum $f^{(i)}$. The associated P-values of the GPF and
$F_{\max}$ tests are displayed in Table \ref{tab.ischemia}. It 
shows that the difference between the normal and ischemia groups was significant, at the 0.05 level, only when we applied the
$F_{\max}$-test. The significance of the $F_{\max}$-test based on the adaptive ECG waveform
indicates that the 12 lead ECG waveform does contain information
about ischemic myocardium. However, this information is masked by
the other physiological facts, like the lead system effect, the HRV
and the respiration. HRV is the phenomenon of non-linear and
non-stationary rhythms of heartbeats, caused by autonomic function,
respiration and other physiological activities, even in resting
conditions \citep{Clifford_Azuaje_McSharry:2006}. Thus, if we
directly perform a Fourier transform to the 12 lead ECG signal, HRV
would distort the result. In addition, lead system effect is
inevitable -- although all leads are placed in fixed locations on
the body surface according to the standard, the cardiac axes vary
among individuals and for each beat. These result in undesirable
signal variations which hamper the visual inspection. Furthermore,
note that the impedance inside the chest wall varies proportional to
the breathing cycle, that is, the impedance increases as we inhale,
and decreases as we exhale. Thus, the respiratory pattern will
further compromise the information. The influence of the respiratory
pattern is clearly evidenced by the insignificance of the $F_{\max}$-test based on 
the adaptive lead signal $\bs_0^{(i,l)}$ as shown in Table
\ref{tab.ischemia}.

The significant results of the $F_{\max}$-test based on the power spectrum $P^{(i)}$ and the cumulative power spectrum $f^{(i)}$ confirm our another hypothesis -- the power spectrum of the manipulated ECG signal is shifted to the low frequency. Indeed, from the middle right panel in Fig. \ref{fig3}, visually we can see the group shift from the high frequency area to the low frequency area. Due to the group shifting phenomenon, the cumulative power spectrum further enhances the significance of the $F_{\textup{max}}$ test, which is shown in Table~\ref{tab.ischemia}.

In conclusion, with the $F_{\textup{max}}$-test, we prove the concept that quantitative analysis of a normalized spectrum from the resting ECG signals contains profound information for screening myocardial ischemia in outpatient settings for those patients with NSR. This finding supports our hypothesis led by the findings in \cite{Mor-Avi_Akselrod:1990,Abboud_Cohen_Selwyn:1987,Pettersson_Pahlm_Carro:2000}. From the clinical viewpoint, the most significant advantage of the approach is that it is safer since it does not entail stress-inducing procedures and it is more convenient, affordable, and cost-effective than the existing methods. We emphasize that for more than a hundred years since the introduction of ECG, physicians have been interpreted ECG in the time domain. In their attempts to diagnose myocardial ischemia using extremely advanced and invasive technologies, researchers have almost completely overlooked more adaptive and accessible approaches to recognize myocardial ischemia. A further large scale clinical study is ongoing and the result will be reported in the near future.

\begin{table}[h]
\begin{center}
\footnotesize\caption{\label{tab.ischemia}  P-values of the  GPF
and $F_{\textup{max}}$ tests for screening stable ischemia. The
P-values of the $F_{\textup{max}}$-test were obtained based on
$10000$ bootstrap replicates. Note that the $F_{\max}$-test is not significant if
we use the non-normalized adaptive lead signal $\bs_0^{(i,l)}$. This
indicates the importance of removing the respiratory effect.}
\begin{tabular}{| c | c | c |}
\hline
&\multicolumn{2}{|c|}{\underline{ANOVA}} \\
 & $F_{\max}$ & GPF   \\ \hline
Adaptive lead signal $\bs_0^{(i,l)}$ & 0.0764 & 0.1317 \\ \hline
Adaptive ECG waveform $\bs^{(i)}$ & 0.0006 & 0.0578   \\ \hline
Power spectrum $P^{(i)}$ & 0.0084 & 0.1566  \\ \hline
Cumulative power spectrum $f^{(i)}$ & 0.0029 & 0.1531   \\
\hline
\end{tabular}
\end{center}
\end{table}

\section{Concluding Remarks}\label{section:conclusion}

In this paper, we proposed and studied the $F_{\max}$-test for the one-way ANOVA problem (\ref{anova.sec1}) for functional data under very general conditions. Via intensive simulation studies, we found that the bootstrapped  null pdf of the $F_{\max}$-test approximates the associated simulated null pdf rather well; in terms of size controlling, the $F_{\max}$-test  outperforms the GPF test of \cite{Zhang_Liang:2013};  and in terms of power, the $F_{\max}$-test generally has higher powers than  the GPF test when the functional data are moderately or highly correlated and the two comparative powers otherwise.

In the preliminary study detailed in Section \ref{section:ischemia:screening}, the $F_{\max}$-test validates the potential of using the adaptive ECG waveform, the power spectrum and the cumulative power spectrum in screening stable ischemia.  A further larger scale study is under way so as to confirm the  results.

Notice that the $F_{\max}$-test is widely applicable. Actually, it is rather easy to extend the $F_{\max}$-test for various hypothesis testing problems for functional data, because for these hypothesis testing problems it is very easy to construct the associated pointwise $F$-test \citep{Ramsay_Silverman:2005,Zhang:2013}. These hypothesis testing problems include two-way ANOVA problems for functional data \citep[ch. 5]{Zhang:2013} and functional linear models with functional responses \citep[ch. 6]{Zhang:2013}, among others. 
\newline\newline
\noindent{\bf Acknowledgments} Zhang's research was supported partially
by the National University of Singapore research grant
R-155-000-108-112 and by the Mathematics Division, National Center
of Theoretical Sciences (Taipei Office). Cheng's research was
supported in part by the National Science Council grants
NSC97-2118-M-002-001-MY3 and NSC101-2118-M-002-001-MY3, and the
Mathematics Division, National Center of Theoretical Sciences
(Taipei Office). Chi-Jen Tseng acknowledges Wei-min Brian Chiu for
the discussion and data collection. Hau-tieng Wu was supported by
the AFOSR grant FA9550-09-1-0643.

\bibliographystyle{biometrika}
\bibliography{fanova1w_ischemia}

\appendix

\setcounter{equation}{0}
\setcounter{section}{0}
\renewcommand{\theequation}{A.\arabic{equation}}
\renewcommand{\thesection}{A.\arabic{section}}
\renewcommand{\thetheorem}{A.\arabic{theorem}}
\renewcommand{\thelemma}{A.\arabic{lemma}}
\renewcommand{\thefigure}{A.\arabic{figure}}
\renewcommand{\thetable}{A.\arabic{table}}

\section{Appendix: Proofs of the Theoretical Results}


To give the proofs of the main results, we first  state the following lemma whose proof is  given in \cite{Zhang_Liang:2013}.

\begin{lemma}\label{lem.sec2} Under Condition A,  as $n\rightarrow \infty$,  we
have \begin{equation}\label{KScv1a.sec2}
 \bz_n(t)\cvgd  GP_k(\bzero, \gamma\bI_k),\;\;
 \sqrt{n}\left\{\hgamma(s,t)-\gamma(s,t)\right\}\cvgd
\GP(0,\varpi),
\end{equation} where  $
 \varpi\left\{(s_1,t_1),(s_2,t_2)\right\}=
       \E\left\{v_{11}(s_{1})v_{11}(t_1)v_{11}(s_2)v_{11}(t_2)\right\}-\gamma(s_1,t_1)\gamma(s_2,t_2).
$ In addition, we have \begin{equation}\label{KScv1b.sec2}
\hgamma(s,t)=\gamma(s,t)+O_{\textup{UP}}\left[n^{-1/2}\right], \end{equation} where
$O_{\textup{UP}}$ means ``bounded in probability uniformly''.
\end{lemma}

We are now ready to outline the proofs of the main results of this paper.

\subsection{Proof of Proposition \ref{pro1.sec2}} Under the given conditions and by Lemma~\ref{lem.sec2},  as $n\rightarrow\infty$, we have $\SSE_n(t)/(n-k)=\hgamma(t,t)\cvgp\gamma(t,t)$ uniformly for
all $t\in\calT$, and  $\bz_n(t)\cvgd \bz(t)\sim
\GP_k(\bzero,\gamma\bI_k)$. Under the null hypothesis, we have $\bmu_n(t)\equiv \bzero$. By Slutsky's theorem,
 as $n\rightarrow\infty$, we have $F_n(t)-R(t)\cvgp 0$ for all $t\in\calT$ where $R(t)=(k-1)^{-1} \bz(t)^T(\bI_k-\bb\bb^T)\bz(t)/\gamma(t,t)$ and
 $\bI_k-\bb\bb^T$ is the limit matrix of $\bI_k-\bb_n\bb_n^T$ as given by  (\ref{blim.sec2}). Since $\calT$ is a finite interval and $F_n(t)$ is continuous over $\calT$, it is also equicontinuous. By Theorem 2.1 in \cite{Newey:1991}, $F_n(t)-R(t) \cvgp 0$ uniformly over $\calT$. Since
 we always have $|\sup_{t\in\calT} F_n(t)-\sup_{t\in\calT} R(t)|\le \sup_{t\in\calT} |F_n(t)-R(t)|$, we have
 $\sup_{t\in\calT} F_n(t)-\sup_{t\in\calT} R(t)\cvgp 0$ which implies that $\sup_{t\in\calT} F_n(t)\cvgd \sup_{t\in\calT} R(t)$. That is
 $F_{\max}\cvgd R_0$ where $R_0= \sup_{t\in\calT} R(t)$.

  Notice that
$\bI_k-\bb\bb^T$ has the singular value decomposition (\ref{bsvd.sec2}). Let
\begin{align}\label{proof_pro1_sec2_eq_bw}
\bw(t)=(\bI_{k-1}, \bzero)\bU^T\bz(t)/\sqrt{\gamma(t,t)}=[w_1(t),w_2(t),\cdots, w_{k-1}(t)]^T.
\end{align}
Then $\bw(t)\sim \GP_{k-1}(\bzero, \gamma_w\bI_{k-1})$
where $\gamma_w(s,t)=\gamma(s,t)/\sqrt{\gamma(s,s)\gamma(t,t)}$. It follows that  $R(t)=(k-1)^{-1}\bw(t)^T\bw(t)=(k-1)^{-1}\sum_{i=1}^{k-1} w_i^2(t)$.   This completes the proof of
Proposition \ref{pro1.sec2}.

\subsection{Proof of Proposition \ref{pro4.sec2}} First of all, notice that given the original $k$ samples (\ref{ksamp.sec1}), the
bootstrapped $k$ samples $v_{ij}^*(t), j=1,2,\cdots, n_i;i=1,2,\cdots, k\iidsim \SP(0,\hgamma)$ where $\hgamma(s,t)$ is the pooled sample covariance function (\ref{hgam.sec2}).  That is to say, the bootstrapped $k$ samples satisfy the null hypothesis (\ref{anova.sec1}).  By Lemma~\ref{lem.sec2} and under Condition A, as $n\rightarrow\infty$,  we have $\hgamma(s,t)\cvgd \gamma(s,t)$ uniformly over $\calT^2$. Applying Proposition \ref{pro1.sec2} leads to the first claim of the proposition and the second claim of the proposition follows immediately. The proposition is then proved.

\subsection{Proof of Proposition \ref{pro2.sec2}} In the proof of Proposition \ref{pro1.sec2}, we have showed that
 under Condition A,   as $n\rightarrow\infty$, we have $\SSE_n(t)/(n-k)\cvgp\gamma(t,t)$ uniformly for
all $t\in\calT$ and  $\bz_n(t)\cvgd \bz(t)\sim
\GP_k(\bzero,\gamma\bI_k)$. Similar to the proof of
Proposition \ref{pro1.sec2}, since $\calT$ is a finite interval and
$F_n(t)$ is equicontinuous over $\calT$, by Slutsky's theorem,
Theorem 2.1 of Newey (1991),  and (\ref{SSRexp2.sec2}), we can show
that as $n\rightarrow\infty$, we have $F_{\max}\cvgd R_1$ with
$R_1=\sup_{t\in\calT}\Big\{(k-1)^{-1}
[\bz(t)+\bd(t)]^T(\bI_k-\bb\bb^T)[\bz(t)+\bd(t)]/\gamma(t,t)\Big\}$
where the idempotent matrix  $\bI_k-\bb\bb^T$  has the singular
value decomposition (\ref{bsvd.sec2}). Let $\bw(t)$ be defined as in
the proof of Proposition ~\ref{pro1.sec2} and let
$\bdelta(t)=(\bI_{k-1},
\bzero)\bU^T\bd(t)/\sqrt{\gamma(t,t)}=[\delta_1(t),\delta_2(t),\cdots,
\delta_{k-1}(t)]^T$. Then $\bw(t)\sim \GP_{k-1}(\bzero,
\gamma_w\bI_{k-1})$ with
$\gamma_w(s,t)=\gamma(s,t)/\sqrt{\gamma(s,s)\gamma(t,t)}$ and
$(\bI_{k-1},\bzero]\bU^T[\bz(t)+\bd(t)]=\bw(t)+\bdelta(t)$.
Therefore, $R_1=\sup_{t\in\calT}
\Big\{(k-1)^{-1}[\bw(t)+\bdelta(t)]^T[\bw(t)+\bdelta(t)]\Big\}=\sup_{t\in\calT}\Big
\{(k-1)^{-1}\sum_{i=1}^{k-1} [w_i(t)+\delta_i(t)]^2  \Big\}$.   This
completes the proof of Proposition \ref{pro2.sec2}.

\subsection{Proof of Proposition \ref{pro3.sec2}}  By (\ref{Tn2Fmax.sec2}), we first have  $P(F_{\max}\ge C_{\alpha}^*)\ge P(T_n\ge  (b-a)C_{\alpha}^*)$. Notice that under Condition A and by Proposition \ref{pro4.sec2},  we have $(b-a)C_{\alpha}^*\cvgd (b-a)C_{\alpha}$ with $C_{\alpha}$ being the upper $100\alpha$ percentile of $R_0$.  Under Condition A and the local alternative (\ref{H1n.sec2}), by the proof of Proposition 3 in \cite{Zhang_Liang:2013}, we have $P(T_n\ge (b-a)C_{\alpha}^*)\rightarrow 1$ as $\delta\rightarrow \infty$.  Proposition \ref{pro3.sec2} is then proved.

\subsection{Proof of Proposition \ref{Fmax_null:discretization}}

By definition, the random vectors
$\bv_{ij,M}=(v_{ij}(t_1),\ldots,v_{ij}(t_M))^T$, $j=1,\ldots,n_i$,
$i=1,2,\cdots, k$, are i.i.d. with a zero mean vector and an
$M\times M$ covariance matrix $\bGamma_M$ whose $(p,q)$ entry is
$\E v_{ij}(t_p)v_{ij}(t_q)=\gamma(t_p,t_q), \, p,q=1,\ldots,M. $
Note that for any finite $M$, taking $n\to\infty$ is exchangeable
with taking maximum over the discretization points $t_1,\ldots,t_M$.
Set $\bz_{n,M}=[\bz_n(t_1)^T,\ldots,\bz_n(t_M)^T]^T$ and
$\bz_M=[\bz(t_1)^T,\ldots,\bz(t_M)^T]^T$ where $\bz_n(t)$ is defined
in (\ref{Zn.sec2}) and $\bz(t)$ is defined in the proof of
Proposition~\ref{pro1.sec2}.  It follows from Lemma \ref{lem.sec2}
that, under Condition A and as $n\rightarrow \infty$, we have
 \begin{equation}\label{KScv1a.discretization}
 \bz_{n,M} \cvgd  \bz_M \sim  N_{kM}(\bzero, \bGamma_M\otimes \bI_k),\;\;
 \sqrt{n}\left\{vec(\hbGamma_M) - vec(\bGamma_M) \right\}\cvgd N_{M^2}(\bzero,\bV_M),
\end{equation} where $\otimes$ is the usual Kronecker product of two matrices,
 $vec(\bA)$ denotes a column vector obtained via stacking all the column
vectors of the matrix  $\bA$ one by one,  and $\bV_M$ is an
$M^2\times M^2$ matrix whose $\left((k_1,l_1),(k_2,l_2)\right)$
entry is
$\E\left[v_{11}(t_{k_1})v_{11}(t_{l_1})v_{11}(t_{k_2})v_{11}(t_{l_2})\right]-\gamma(t_{k_1},t_{l_1})\gamma(t_{k_2},t_{l_2})$,
and $v_{11}(t)$ is the subject-effect function of the first subject
of the first group. In addition, we have
\begin{equation}\label{KScv1a2.discretization} \hbGamma_M =
\bGamma_M+O_{\textup{UP}}\big[n^{-1/2}\big]. \end{equation} By
(\ref{KScv1a.discretization}) we have
$ \SSR_{n}(t_l)/(k-1)\cvgd \bz(t_l)^T(\bI_k-\bb\bb^T)\bz(t_l)/(k-1).
$
Under the given conditions and by (\ref{KScv1a.discretization}) and (\ref{KScv1a2.discretization}), as $n\rightarrow\infty$, we have
$ \SSE_{n}(t_l)/(n-k)=\hgamma(t_l,t_l)\cvgp\gamma(t_l,t_l) $
 for all $l=1,\ldots,M$.
Under the null hypothesis and by  Slusky's Theorem, as
$n\rightarrow\infty$, we have
$F_{\max,M}=\displaystyle{\max_{l=1,\ldots,M}}
\{\SSR_{n}(t_l)/(k-1)\}\{\SSE_{n}(t_l)/(n-k)\} \cvgd R_{0,M}$ where
$R_{0,M}$ is defined by
\[
R_{0,M}=\max_{l=1,\ldots,M}\Big[\big\{(k-1)\gamma(t_l,t_l)\big\}^{-1}\bz(t_l)^T(\bI_k-\bb\bb^T)\bz(t_l)\Big],
\]
and $\bI_k-\bb\bb^T$ is the limit matrix of $\bI_k-\bb_n\bb_n^T$;
see  (\ref{blim.sec2}). Note that under the null hypothesis, we have
$ \bmu(t_l)^T(\bI_k-\bb\bb^T)\bmu(t_l) \equiv \bzero, l=1,\ldots,
M$. For $l=1,\ldots,M$, set
\begin{align}\label{discretization:proof:definition:wM}
\bw_M(l) =(\bI_{k-1}, \bzero)\bU^T\bz_M(l) \big/
\sqrt{\gamma(t_l,t_l)}=[w_{1,M}(l),w_{2,M}(l),\cdots,w_{k-1,M}(l)]^T,
\end{align}
where $\bU$ comes from the singular value decomposition
(\ref{bsvd.sec2}) of $\bI_k-\bb\bb^T$. Then we have  $
R_{0,M}=\max_{l=1,\ldots,M}
\big\{(k-1)^{-1}\bw_M(l)^T\bw_M(l)\big\}=\max_{l=1,\ldots,M}\big
\{(k-1)^{-1}\sum_{i=1}^{k-1} w_{i,M}^2(l)  \big\}.$ Let $\bw_{i,M}=(w_{i}(t_1),\ldots,w_{i}(t_M))^T, i=1,\cdots, k-1$. Then
$\bw_{1,M},\ldots,\bw_{k-1,M}\iidsim N_M(\bzero,\bGamma_{w,M})$
where the $(p,q)$ entry of $\bGamma_{w,M}$ is $\gamma_w(t_p,t_q),
p,q=1,\ldots,M$.
This completes the proof of $F_{\max,M}\cvgd R_{0,M}$ as $n\rightarrow\infty$. 

Next, by Condition A5 on $\gamma(s,t)$ 
and a direct calculation by the definition of H\"older's continuity,
the covariance function $\gamma_w(s,t)$ of the Gaussian process
$\omega_i(t)$ defined in (\ref{proof_pro1_sec2_eq_bw}) is in
$C^{\beta}(\calT\times\calT)$. Thus, Kolmogorov's Theorem
\cite[Theorem 18.19]{Koralov_Sinai:2007} says that, for all
$i=1,\ldots,k-1$, there exists a continuous modification $W_i(t),
t\in\calT$
of $\omega_i(t), t\in\calT$ 
and an event subspace $\Omega(M)$ with $P(\Omega(M))\to 1$ as $M\to \infty$, so that for all events in $\Omega(M)$, $W_i(t)$ is H\"older continuous
with exponent $0<\tilde{\beta}< \beta$ and H\"older's modulus $c_{\tilde{\beta}}=2/(1-2^{-\tilde{\beta}})$. 
Thus, from now on we work with $\{W_i(t)\}_{i=1}^{k-1}$ instead of
$\{\omega_i(t)\}_{i=1}^{k-1}$ and use the same notation for the
versions of $R_{0}$ and $R_{0,M}$ based on the continuous
modifications:
\begin{align*}
R_{0} \dequ \sup_{t\in\mathcal{T}}\Big\{ (k-1)^{-1} \sum_{i=1}^{k-1}
W_i^2(t)\Big\},\quad R_{0,M} \dequ \max_{l=1,\ldots,M}\Big\{
(k-1)^{-1} \sum_{i=1}^{k-1} W_{i,M}^2(l)\Big\},
\end{align*}
where $W_{i,M}(l)=W_i(t_l), l=1,\ldots, M; i=1,\ldots,k-1$. Clearly
$(k-1)^{-1}\sum_{i=1}^{k-1} W_i^2(t)$ is also H\"older continuous
with exponent $\tilde{\beta}$ and H\"older's modulus $\tilde{c}=
c_{\tilde{\beta}}^2$.
Since $\mathcal{T}$ is compact, the supremum of $\sum_{i=1}^{k-1} W_i^2(t)$ is achieved at some 
$t'\in[t_{l'-1},t_{l'+1}]$. By H\"older's continuity of the process, we have
\[
\Big|(k-1)^{-1}\sum_{i=1}^{k-1}
W_{i,M}^2(l')-\sup_{t\in[t_{l'-1},t_{l'+1}]}
(k-1)^{-1}\sum_{i=1}^{k-1} W_i^2(t) \Big|\leq\tilde{c}
\tau_M^{\tilde{\beta}}.
\]
Suppose the maximum of $\sum_{i=1}^{k-1} W_{i,M}^2(l)$ is achieved
at $l''$ instead of $l'$, we must have
\[
 (k-1)^{-1} \sum_{i=1}^{k-1} W_{i,M}^2(l') \leq (k-1)^{-1} \sum_{i=1}^{k-1} W_{i,M}^2(l'') \leq (k-1)^{-1} \sum_{i=1}^{k-1} W_i^2(t'),
\]
where the last inequality holds by definition.
Thus, for all events in $\Omega(M)$, we have
\begin{equation}\label{discretization:diff}
\big|R_0- R_{0,M}\big|\leq \tilde{c} \tau_M^{\tilde{\beta}}.
\end{equation}

\subsection{Proof of Proposition \ref{discretization:NPB}}

Note that the $k$ bootstrapped samples  $\{\bv_{ij,M}^*\}_{j=1}^{n_i},i=1,\cdots, k$ are i.i.d. with mean vector $\bzero$ and covariance
matrix $\hbGamma_M$ whose $(p,q)$ entry is $\hgamma(t_p,t_q), p,q=1,\ldots,M$. That is to say, the bootstrapped $k$ samples satisfy the
discretized null hypothesis (\ref{discrete:H0}). 
By (\ref{KScv1a.discretization}) and (\ref{KScv1a2.discretization}),
and under Condition A, as $n\rightarrow\infty$,  we have $\hbGamma_M
\cvgd \bGamma_M$ uniformly. Applying the same arguments as in the
proof of $F_{\max,M}\cvgd R_{0,M}$ to the $k$ bootstrapped samples
$\{\bv_{ij,M}^*\}_{j=1}^{n_i},i=1,\cdots, k$, leads to the claim
that $F_{\max,M}^*\cvgd R_{0,M}$ as $n\rightarrow\infty$.

\subsection{Proof of Proposition \ref{discretization:pro2.sec2}}

As in the proof of Proposition \ref{Fmax_null:discretization}, under
Condition A, $\SSE_{n}(t_l)/(n-k)\cvgp\gamma(t_l,t_l)$ and
$\bz_{n,M} \cvgd \bz_M\sim N_{kM}(\bzero,\bGamma_M\otimes \bI_k)$ as
$n\rightarrow\infty$. Under $H_1$, we have $ \mu_i(t_l) =\mu_0(t_l)
+n_i^{-1/2} d_i(t_l) , l=1,\ldots, M;\; i=1,2,\cdots, k.$ It follows
that  $ \bmu(t_l) =\mu_0(t_l)\mathbf{1}_k+\diag(n_1^{-1/2},\ldots,
n_k^{-1/2}) \bd(t_l),\, l=1,\ldots,M$ where
$\bmu(t)=[\mu_1(t),\ldots,\mu_k(t)]^T$ and
$\bd(t)=[d_1(t),\ldots,d_k(t)]^T$ as defined before. Then we have
$\bmu_n(t_l)=\diag(n_1^{1/2},\ldots,n_k^{1/2})\bmu(t_l)=\mu_0(t_l)\bb_n+\bd(t_l).$
Since $(\bI_k-\bb_n\bb_n^T/n)\bb_n=\bzero$, under $H_1$, we have
$ \SSR_{n}(t_l) =\big[\bz_n(t_l)
+\bd(t_l)\big]^T(\bI_k-\bb_n\bb_n^T/n)\big[\bz_n(t_l)+\bd(t_l)\big].
$
Hence, as $n\rightarrow\infty$, we have $F_{\max,M}\cvgd R_{1,M}$.
Let $\bw_M(l)$ be as defined in
(\ref{discretization:proof:definition:wM}) and let $
\bdelta_M(l)=\bdelta(t_l)=(\bI_{k-1},
\bzero)\bU^T\bd(t_l)/\sqrt{\gamma(t_l,t_l)}=[\delta_{1,M},\delta_{2,M},\cdots,
\delta_{k-1,M}]^T $ where $\delta_{i,M}(l)=\delta_i(t_l),
l=1,\ldots,M$. Thus,
$(\bI_{k-1},\bzero]\bU^T[\bz(t_l)+\bd(t_l)]/\sqrt{\gamma(t_l,t_l)}=\bw_M(l)+\bdelta_M(l)$.
Therefore, we have 
\begin{eqnarray*}
 R_{1,M}&=&\max_{l=1,\ldots,M}
\Big\{(k-1)^{-1}[\bw_M(l)+\bdelta_M(l)]^T[\bw_M(l)+\bdelta_M(l)]\Big\}\\
&=&\max_{l=1,\ldots,M}\Big \{(k-1)^{-1}\sum_{i=1}^{k-1}
[w_{i,M}(l)+\delta_{i,M}(l)]^2  \Big\}. 
\end{eqnarray*}

Next, using the same arguments in the proof of (\ref{discretization:diff}), 
we can show that there exists an event space $\Omega(M)$ with $P(\Omega(M))\rightarrow 1$ such that conditional on $\Omega(M)$,
for all $i=1,\ldots,k-1$, we have a continuous modification $W_i(t)$ of $\omega_i(t)$ defined in (\ref{proof_pro1_sec2_eq_bw}) so that $W_i(t)$
is H\"older continuous with exponent $0<\tilde{\beta}< \beta/2$ 
and H\"older's modulus $c_{\tilde{\beta}}=2/(1-2^{-\tilde{\beta}})$. We use the same notation for the versions of $R_{1}$ and $R_{1,M}$
 based on the continuous modifications:
\begin{align}
&R_{1} \dequ  \sup_{t\in\calT}\Big \{(k-1)^{-1}\sum_{i=1}^{k-1} [W_i(t)+\delta_i(t)]^2  \Big\},\,\, 
R_{1,M} \dequ \max_{l=1,\ldots,M}\Big\{ (k-1)^{-1} \sum_{i=1}^{k-1}
[W_{i,M}(l)+\delta_{i,M}(l)]^2\Big\},\nonumber
\end{align}
where $W_{i,M}(l)=W_i(t_l), l=1,\ldots,M$.
Clearly, under the assumption on $d_i(t)$, we know $(k-1)^{-1}\sum_{i=1}^{k-1} [W_i(t)+\delta_i(t)]^2$ is
also H\"older continuous with exponent $\tilde{\beta}$ and H\"older's modulus $\tilde{c} = 2(c_{\tilde{\beta}}^2+c^2_{\delta})$,
where $c_{\delta}$ is the maximum of the H\"older's modulus of $d_i(t)$, $i=1,\ldots,k$. Again,
by the same arguments as those leading to (\ref{discretization:diff}), 
for all events in $\Omega(M)$ we have
$
\big|R_1- R_{1,M}\big|\leq \tilde{c} \tau_M^{\tilde{\beta}}.
$

\end{document}